\DeclareSymbolFont{AMSb}{U}{msb}{m}{n}
\documentclass[11pt,a4paper,leqno,noamsfonts]{amsart}

   \makeatletter
   \renewcommand\@biblabel[1]{#1.}
   \makeatother
   \usepackage{multirow}

\usepackage{graphicx}
\usepackage{subcaption}

\usepackage{fancyvrb}

\usepackage{mathtools}

\usepackage[english]{babel}
\usepackage[dvipsnames]{xcolor}
\usepackage{graphicx,pifont,soul,physics} 
\usepackage[utopia]{mathdesign}
\usepackage[a4paper,top=3cm,bottom=3cm,
           left=3cm,right=3cm,marginparwidth=60pt]{geometry}

\usepackage{fancyhdr}
\usepackage{float}

\usepackage[utf8]{inputenc}
\usepackage{braket,caption,comment,mathtools,stmaryrd}
\usepackage[usestackEOL]{stackengine}
\usepackage{multirow,booktabs,microtype,relsize}
\usepackage[colorlinks,bookmarks]{hyperref} %
      \hypersetup{colorlinks,%
           citecolor=britishracinggreen,%
            filecolor=black,%
            linkcolor=cobalt,%
            urlcolor=cornellred}
      \numberwithin{equation}{section}
\usepackage[capitalise]{cleveref}
\usepackage[colorinlistoftodos]{todonotes}
\definecolor{antiquewhite}{rgb}{0.98, 0.92, 0.84}
\definecolor{buff}{rgb}{0.94, 0.86, 0.51}
\definecolor{palecopper}{rgb}{0.85, 0.54, 0.4}
\definecolor{fluorescentyellow}{rgb}{0.8, 1.0, 0.0}
\definecolor{bole}{rgb}{0.47, 0.27, 0.23}

\usepackage{amsmath,float}
\usetikzlibrary{decorations.markings}

\usepackage[vcentermath]{youngtab}
\input xy
\xyoption{all}

\usepackage{pgfkeys}
\usepackage{pgfopts}
\usepackage{ytableau}

\definecolor{cornellred}{rgb}{0.7, 0.11, 0.11}
\definecolor{britishracinggreen}{rgb}{0.0, 0.26, 0.15}
\definecolor{cobalt}{rgb}{0.0, 0.28, 0.67}
\DeclareSymbolFont{usualmathcal}{OMS}{cmsy}{m}{n}
\DeclareSymbolFontAlphabet{\mathcal}{usualmathcal}

\newcommand{\BA}{{\mathbb{A}}}

\newcommand{\BC}{{\mathbb{C}}}

\newcommand{\BR}{{\mathbb{R}}}

\newcommand{\BZ}{{\mathbb{Z}}}

\newcommand{\CN}{{\mathcal N}}
\newcommand{\CO}{{\mathcal O}}

\newcommand{\CS}{{\mathcal S}}

\newcommand{\mult}{\mathrm{mult}}

\newcommand{\vol}{{\mathrm{vol}}}
\newcommand{\area}{{\mathrm{area}}}

\DeclareMathOperator{\length}{length}

\DeclareFontFamily{OT1}{rsfs}{}
\DeclareFontShape{OT1}{rsfs}{n}{it}{<-> rsfs10}{}
\DeclareMathAlphabet{\curly}{OT1}{rsfs}{n}{it}
\renewcommand\hom{\mathscr{H}\kern-0.3em\mathit{om}}

\DeclareMathOperator{\lHom}{\mathscr{H}\kern-0.3em\mathit{om}}
\DeclareMathOperator{\RRlHom}{\mathbf{R}\kern-0.025em\mathscr{H}\kern-0.3em\mathit{om}}

\DeclareMathOperator{\lExt}{{\mathscr{E}\kern-0.2em\mathit{xt}}}

\newcommand\Spec{\operatorname{Spec}}

\newcommand\Proj{\operatorname{Proj}}
\newcommand\Bl{\operatorname{Bl}}



\usepackage[all]{xy}
\usepackage{tikz}
\usepackage{tikz-cd}
\usepackage{adjustbox}
\usepackage{rotating}
\usepackage{comment}

\usetikzlibrary{matrix,shapes,intersections,arrows,decorations.pathmorphing}
\tikzset{commutative diagrams/arrow style=math font}
\tikzset{commutative diagrams/.cd,
mysymbol/.style={start anchor=center,end anchor=center,draw=none}}

\tikzset{
shift up/.style={
to path={([yshift=#1]\tikztostart.east) -- ([yshift=#1]\tikztotarget.west) \tikztonodes}
}
}

\theoremstyle{definition}

\newtheorem*{lemma*}{Lemma}
\newtheorem*{theorem*}{Theorem}
\newtheorem*{example*}{Example}
\newtheorem*{fact*}{Fact}
\newtheorem*{notation*}{Notation}
\newtheorem*{definition*}{Definition}
\newtheorem*{prop*}{Proposition}
\newtheorem*{remark*}{Remark}
\newtheorem*{corollary*}{Corollary}

\newtheorem*{conventions*}{Conventions}
\newtheorem*{question*}{Question}
\newtheorem*{convention*}{Convention}

\newtheorem{definition}{Definition}[section]

\newtheorem{example}[definition]{Example}

\newtheorem{remark}[definition]{Remark}

\newtheoremstyle{thm} 
        {3mm}
        {3mm}
        {\slshape}
        {0mm}
        {\bfseries}
        {.}
        {1mm}
        {}
\theoremstyle{thm}
\newtheorem{theorem}[definition]{Theorem}
\newtheorem{corollary}[definition]{Corollary}
\newtheorem{lemma}[definition]{Lemma}

\newtheorem{prop}[definition]{Proposition}
\newtheorem{thm}{Theorem}

\newtheoremstyle{ex} 
        {3mm}
        {3mm}
        {}
        {0mm}
        {\scshape}
        {.}
        {1mm}
        {}
\theoremstyle{ex}

\newtheoremstyle{sol} 
        {3mm}
        {3mm}
        {}
        {0mm}
        {\scshape}
        {.}
        {1mm}
        {}
\theoremstyle{sol}

\usepackage{tikz}
\usepackage{xparse}

\newcount\tableauRow
\newcount\tableauCol

\newenvironment{Tableau}[1]{%
  \tikzpicture[scale=0.5,draw/.append style={thick,black},
                      baseline=(current bounding box.center)]
    \tableauRow=-1.5
    \foreach \Row in {#1} {
       \tableauCol=0.5
       \foreach\k in \Row {
         \draw[thin](\the\tableauCol,\the\tableauRow)rectangle++(1,1);
         \draw[thin](\the\tableauCol,\the\tableauRow)+(0.5,0.5)node{$\k$};
         \global\advance\tableauCol by 1
       }
       \global\advance\tableauRow by -1
    }
}{\endtikzpicture}

\newtheorem*{Acknowledgments*}{Acknowledgments}
\DeclareMathAlphabet\BCal{OMS}{cmsy}{b}{n}

\address{Centre for Mathematical Sciences, University of Cambridge, Wilberforce Road, CB3 0WA, Cambridge, United Kingdom}
\email{fr414@cam.ac.uk}

\title[A quotient formula for the Behrend function]{A quotient formula for the Behrend function}
\author{Fatemeh Rezaee}
\date{}

\begin{document}
\maketitle

\begin{abstract}
The Behrend function is essential in enumerative geometry.  However, computing it, even in zero dimension, has been a highly subtle and mysterious matter since its introduction in the 2000s. 
We present a quotient formula for the Behrend function of any closed $r$-dimensional affine scheme $X^r \subset \mathbb{A}^N$, evaluated at a singular point, where $0\leq r \leq N-1$.
In particular, if the scheme is monomial, we associate parallelepipeds to it and interpret the Behrend function as simple volume information about them.  The results here, in particular, recover the recent results of Graffeo-Ricolfi for the zero-dimensional case in 2D in \cite{Graffeo-Ricolfi}, and also the more recent result of Polini-Sammartan-Ulrich on the zero-dimensional case in \cite{PSU26}.

\end{abstract}

{\hypersetup{linkcolor=black}
\tableofcontents}

\section{Introduction}
The Behrend function (introduced by Kai Behrend in \cite{Beh}) plays a central role in enumerative geometry; in particular, the Behrend function of zero-dimensional subschemes is crucial in counting invariants in enumerative theories with symmetric perfect obstruction theories (\cite{BFinc,BFHilb}). However, this function is infamous for having unpredictable behaviour; hence, it has been very hard to compute it, even in low dimensions. We provide a practical and elementary way to calculate the Behrend function of any affine scheme at a singular point; the case becomes more interesting when it comes to monomial (affine) schemes\footnote{A \emph{monomial affine scheme} is given by $\Spec(\BC[x_1,\ldots,x_N]/I)$, where $I$ is a monomial ideal.}, which come with a combinatorial intuition. 

The case of zero-dimensional subschemes in $\BA^2$ was first considered in \cite{Graffeo-Ricolfi} by Graffeo and Ricolfi using toric geometry, which covers many interesting non-trivial classes and confirms the difficulties of the Behrend function calculation. Also, more recently, the case of zero-dimensional subschemes in $\BA^N$ was  considered in \cite{PSU26}.  For more information on the applications of the Behrend function in enumerative geometry, we also refer to \cite{Graffeo-Ricolfi}. For the application of the Behrend function in Donaldson-Thomas theory, we refer to \cite{JoyceSong}.

Before we state the main results, we give some quick key definitions (more details will be in Section \ref{subsection:Behrend}.

\begin{definition}[Normal cone and special fibre]\label{Def:normalCone}
Let $I$ be an $r$-dimensional\footnote{By \emph{dimension} of an ideal, we mean the Krull dimension of its quotient ring.} ideal in $A=\BC[x_1,\ldots,x_N]$ and let $X:=\Spec(A/I)$ be singular at the origin $\mathbf{o}\in X$. We define the geometric \textit{normal cone of $I$} to be $\CN(I):=\Spec(\mathrm{gr}(I))=\Spec(\bigoplus_{i=0}^{\infty}\frac{I^i}{I^{i+1}})$. Note that $\Proj(\bigoplus_{i=0}^{\infty}\frac{\mathfrak{m}^i}{\mathfrak{m}^{i+1}})$ defines the exceptional locus of $\Bl_{\mathbf{o}}(X)$.

We, respectively, denote by $\CN_i$ and $\CN_{i,j}$, the maximum-dimensional irreducible components of  ${\CN(I)}$, for $I$, respectively, zero-dimensional or positive-dimensional. We denote by $\mult(I)$ the normalised leading coefficient of the Hilbert polynomial of $A/I$.

Also,  we define the \emph{special fibre of $I$}, $\CS(I)$, to be  the reduction \footnote{For a quotient ring $A/J$, its \emph{reduction} is given by $A/\sqrt{J}$. We use the term reduction for either ideals/quotient rings or their corresponding affine or (quasi-)projective schemes} of $\CN(I)$\footnote{As in the proof of Theorem \ref{v as quotien a of mult (reducedness condition)},  we can write the normal cone as $\Spec(\frac{\BC[x_1,\ldots,x_N,w_1,\ldots,w_l]}{L})$, for an ideal $L$}.

We use the term \emph{degree} for either ideal/quotient rings or their corresponding affine or (quasi-)projective schemes\footnote{The degree of an ideal $I$ or its quotient ring $A/I$ is the factorial of the (Krull) dimension of $A/I$, times the leading coefficient of the corresponding Hilbert polynomial. See Definition \ref{def: degree}.}.
\end{definition}

\begin{definition}[Behrend function, see definitions \ref{Def: Behrend}, \ref{def: distinguished}, \ref{Def?:Behrend2}] Let $X=\Spec(\BC[_1,\ldots,x_N]/I)$ be an affine scheme. We define the Behrend function on $X$ as  $v_X:=\mathrm{Eu}(\mathfrak{c}_X)$, where 
\begin{align*}
\mathfrak{c}_{X}=\mathfrak{c}_{X/{\BA^N}}:=\sum_{\CN_{i,j}\subset\CN(I)}(-1)^{\dim\pi(\CN_{i,j})}l(\CN_{i,j})\big[\pi(\CN_{i,j})\big],
     \end{align*}
     where $l(\CN_{i,j})$ is the length of $\CN(I)$ at the generic point of $\CN_{i,j}$,
     and $\mathrm{Eu}$ is a function from the integral algebraic cycles on $X$ to constructible $\BZ$-valued
functions on $X$, and $\pi: \CN(I)\to X$ is the \emph{canonical projection} map.

\end{definition}

\subsection{Main results} We prove the following main results.

{\bf Setup.} Let $\CN(I)$ be the geometric normal cone of an $r$-dimensional 
ideal $I$ in the polynomial ring $A=\BC[x_1,...,x_N]$, as above. As in the proof of Theorem \ref{v as quotien a of mult (reducedness condition)}, we see that we can write the normal cone as $\Spec(\frac{\BC[x_1,\ldots,x_N,w_1,\ldots,w_l]}{L})$, for an ideal $L$. We can write $L=\cap_{i,j} Q_{i,j}$, where $\{Q_{i,j}\}$ is its $p_{i,j}$-primary decomposition.  Let $\mult_{i,j}$ be the degree of  $Q_{i,j}$.  Let $D_{i,j}$ be the degree of $p_{i,j}$, for $0\leq j \leq r$, and $1\leq i \leq n_j(I)$, where $n_i(I)$ is the number of the irreducible components of $\CN(I)$ of the maximal dimension, dominating an $i$-dimensional component of $\Spec(A/I)$. In addition, we set $n(I)=\sum_{i=0}^{r}n_i(I)$, which is the number of irreducible components of the maximal dimension of $\CN(I)$.

In this article, for simplicity, we assume that the origin $\mathbf{o}\in X$ is a singular point of $X=\Spec(A/I)$, but the results can be considered for any singular point. We are interested in evaluating $v_X$ at the origin $\mathbf{o}$. When $\dim X=0$, we can simply use $v_X$ for $v_X(\mathbf{o})$. 

\begin{thm}[Quotient formula; see Theorem \ref{v as quotien a of mult (reducedness condition)}] \label{generalConjecture} 
For an $r$-dimensional  
ideal $I$ in $\text{\:}A=\BC[x_1,...,x_N]$, suppose  $X=\Spec(A/I)$
has a singularity at the origin $\mathbf{o}$. Then we have

\begin{align}
v_X(\mathbf{o})=&\sum_{j=0}^{r}\sum_{i=1}^{n_j(I)}(-1)^j\frac{\mult_{i,j}}{D_{i,j}}\mathrm{Eu}(\big[\pi(\CN_{i,j})\big]),
\end{align}
where $\CN_{i,j}$ are the maximal-dimensional irreducible components of $\CN(I)$, and $\pi:\CN(I)\to X$ is the canonical projection, and $\:{}\mathrm{Eu}$ is the local Euler obstruction.

\end{thm}

An interesting combinatorial interpretation of this quotient formula for monomial ideals is described below. For a monomial ideal $I$ in $N$ variables, we associate the \emph{convex hull}, $\mathrm{Conv}(I)$, as the convex hull of the exponents of the monomial terms of $I$ in $\BR^N$.

\begin{thm}[A combinatorial interpretation for monomial ideals; see Theorem \ref{combinatorialinterpretation}]  \label{gcdMainThm} Let $0\leq j \leq r$ and  $1\leq i \leq n_j(I)$.
 Let $I$ be an $r$-dimensional monomial ideal in $A=\BC[x_1,\ldots,x_N]$, and let $\mathbf{o}$ be the origin as a singular point of  $\text{\:}X=\Spec(A/I)$. Then 
\begin{align}
v_{X}(\mathbf{o})=&\sum_{j=0}^{r}\sum_{i=1}^{n_j(I)}(-1)^j\gcd_k(\vol(\text{parallelepipeds over  ${T^k_{i,j}}))$}\mathrm{Eu}(\big[\pi(\CN_{i,j})\big]),
\end{align}
 where $\{T^k_{i,j}\}_k$ is a triangulation of $P_{i,j}$, where the latter is the projection of $\text{\:}\mathrm{Conv}(I)$ onto the $(N-j)$-coordinate plane for every fixed $j$. Also, $\CN_{i,j}$ are the maximal-dimensional irreducible components of $\CN(I)$, and $\pi:\CN(I)\to X$ is the canonical projection, and $\:{}\mathrm{Eu}$ is the local Euler obstruction. Also, $P_{i,j}$ corresponds to an irreducible component of the normal cone with maximal dimension, dominating a $j$-coordinate plane component of $\text{\:}\Spec(A/I)$ in $\BA^N$.   The parallelepiped over $T^k_{i,j}$ is a parallelepiped produced by the vectors from the origin on $\BR^N$ to the vertices of $T^k_{i,j}$ (see Figure \ref{Fig:parallelepiped} and Figure \ref{Fig:firstEx}).
  \begin{figure}[H]
 \centering
    \includegraphics[width=6.5cm]{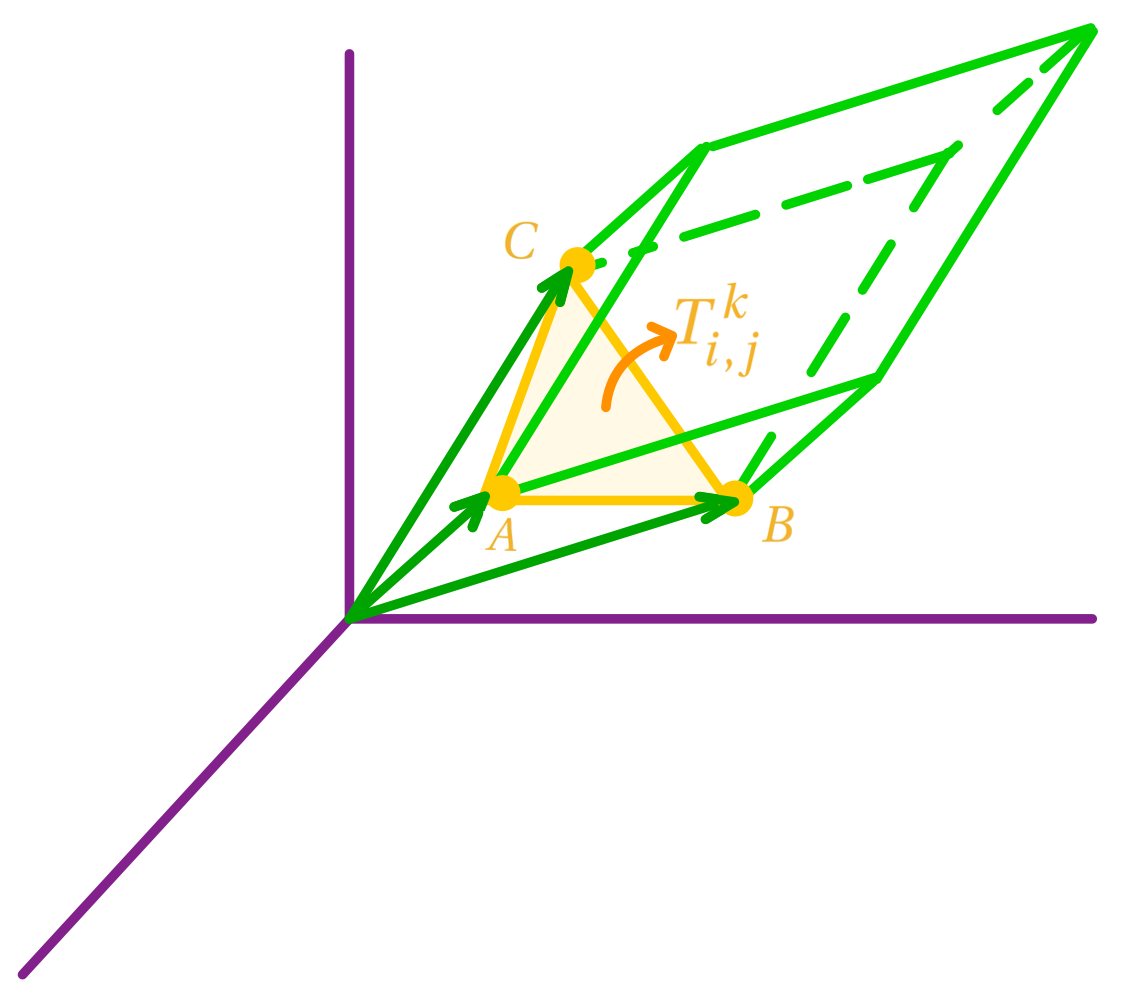}
    \caption{Let the simplex $T^k_{i,j}$ correspond to exponents associated to vertices $A,B$ and $C$. The parallelepiped over the simplex $ABC$ is the parallelepiped generated over the vectors from the origin of $\BR^N$ to the vertices  $A, B$ and $C$.}
\label{Fig:parallelepiped}
\end{figure}
\end{thm}

\begin{remark}\label{coordinatePlanes}
    Note that since in the statement of Theorem \ref{gcdMainThm}, we assumed that $X=\Spec(A/I)$ has the origin as its singularity, the components of $X$ correspond to $i$-coordinate planes (i.e., $i$-dimensional coordinate planes). Therefore, the last sentence of the theorem makes sense.
\end{remark}

\begin{remark}\label{rem: EulerOne}
    Note that in many interesting cases, for instance when  $\pi(\CN_{i,j})$ is smooth at the singularity point $\mathbf{o}\in X$, or it is a single point, by Lemma \ref{PropertiesEuler}, we have both formulas in Theorem \ref{generalConjecture} and Theorem \ref{gcdMainThm} without the contribution the local Euler obstruction. Also, when $r=0$, we always have $\mathrm{Eu}=1$.
\end{remark}

To clarify How Theorem  \ref{generalConjecture} and Theorem \ref{gcdMainThm} work, we give an explicit example.

\begin{example} Let $I=(x^3,y^4,z^6,xy,xz,y^2z^2,yz^4,y^3z)$ be a zero-dimensional ideal and $X=\Spec(A/I)$, which contains the origin as its singularity. We can see that  there are three irreducible components of the normal cone $\CN(I)$ of
the maximal dimension. From Theorem  \ref{generalConjecture}, we have
\begin{align*}
    v_X=&\frac{3}{1}+\frac{12}{3}+\frac{12}{2}=3+4+6=13.
\end{align*}
Since $I$ is a monomial ideal, alternatively from Theorem \ref{gcdMainThm}, we have
\begin{align*}
v_{X}=\gcd(3)+\gcd(4,4,4)+\gcd(6,6)=3+4+6=13.
\end{align*}
 See Figure \ref{Fig:firstEx}.
     \begin{figure}[H]
 \centering
    \includegraphics[width=6.5cm]{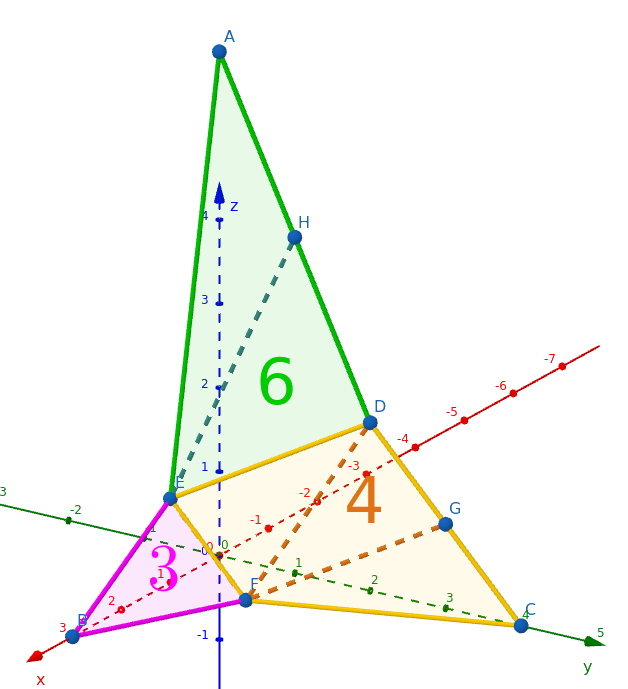}
    \caption{Irreducible components of the normal cone of maximal dimension are pictured with different colours. We have also picture the triangulation corresponding to each irreducible components. The numbers attached to each component are the corresponding gcd's.}
\label{Fig:firstEx}
\end{figure}
\end{example}

An even simpler interpretation of the Behrend function for monomial schemes is the following.

\begin{thm}[see Theorem \ref{v=|d|}]\label{ThmC} Let $I$ be a monomial ideal in $A=\BC[x_1,x_2,\ldots,x_N]$ that is singular at the origin $\mathbf{o}$\footnote{For an ideal $I$ in $A$, we interchangeably say either $\Spec(A/I)$ is singular at a point, or $I$ is.}, and let $X=\Spec(A/I)$. Then, for any polygon $P_{i,j}$ in the lower boundary of the convex hull of $I$, 
we have
\begin{align*}
    v_{X,i,j}(\mathbf{o})=&|d_{i,j}|,
\end{align*}
where $v_{X,i,j}(\mathbf{o})$ is the contribution of the Behrend function evaluated at the origin (considering with a suitable sign), and $d$ is the constant term of the primitive equation $a_1x_1+a_2x_2+\ldots+a_Nx_N+d_{i,j}=0$ of the plane defined by $P_{i,j}$. Then $v_X(\mathbf{o})$ can be computed as a signed sum of $v_{X,i,j}(\mathbf{o})$s, weighted with the corresponding $\mathrm{Eu}(\big[\pi(\CN_{i,j})\big])$.
\end{thm}

To prove Theorem \ref{gcdMainThm}, we need the following theorem, which will be proved inductively.
\begin{thm}[Normal cone/convex hull correspondence; see Theorem \ref{thm:Proof normalconConvexhulldecomposition}] \label{thm: normalconConvexhulldecomposition}
    Let $I$ be an $r$-dimensional monomial ideal. The maximal dimensional irreducible components of the normal cone of $I$ correspond to the polytope decomposition of the lower boundary of the convex hull of $I$ in each $(N-j)$-coordinate plane for $\text{\:}0\leq j \leq r$.
\end{thm}

From Theorem \ref{gcdMainThm} and the \textit{minimisation} procedure of a monomial ideal (we denote the minimisation of $I$ by $I_{\min}$), we can characterise zero-dimensional monomial ideals with a maximal Behrend number, i.e., when the Behrend number is equal to the multiplicity:
  \begin{corollary} [Characterisation of zero-dimensional monomial ideals with maximal Behrend number; see Proposition \ref{maxB}]
  A zero-dimensional monomial ideal $I$ in $\BC[x_1,\ldots,x_N]$ has maximal Behrend number if and only if the $(N-1)$-boundary of the minimisation $I_{\min}$ only consists of $(N-1)$-simplices and the removed monomials in the process of minimisation lie strictly above the lower boundary of the convex hull of $I$.
  \end{corollary}
  Note that this corollary can be generalised to higher dimensions, but we skip the statement.

Our procedure shares the spirit of the following combinatorial result, which will be proved in the appendix. More precisely, we have a generalised version of the celebrated Pick's Theorem that relaxes the condition that the vertices be integral for certain types of plane polygons.

\begin{thm} [Generalised Pick's formula; see Theorem \ref{Pick's}]
    Let $P$ be a plane convex polygon with some non-lattice vertices so that it can be integrally and generically embedded in $\BR^3$, i.e., there is an embedding of $P$ into $\BR^3$ with all vertices lattice points so that it is not perpendicular to any of the three coordinate planes.  Let $\widetilde{I(P)}$ and $\widetilde{B(P)}$ be the number of internal lattice points and the number of lattice points on the boundary of the embedded polygon in $\BR^3$, respectively.  Then
    \begin{align*}
        \area(P)=&\left(\frac{\widetilde{B(P)}}{2}+\widetilde{I(P)}-1\right)\left(\sqrt{a^2+b^2+c^2}\right),
    \end{align*}
    where $ax+by+cz+d=0$ is the equation of the plane defined by the embedded polygon in $\BR^3$, $(a,b,c)$ is primitive, and $d\neq 0$.
\end{thm}

\begin{remark}
    A generalised version for $d=0$ is given in \cite{Si}.
\end{remark}

Richard Thomas raised the following question, but we do not know the answer, and it could be an interesting future direction.
\begin{question*}
    What is the enumerative interpretation of the Behrend function of ideals that do not define a critical locus of a function?
\end{question*}

\subsection*{Previous work} In \cite{Graffeo-Ricolfi}, the authors used toric geometry to compute the Behrend number of fat points in $\BA^2$. Our approach recovers their formulas and extends them to higher dimensions using combinatorial methods. Regarding the constancy of the Behrend function on the Hilbert scheme of points in 3D, we refer to \cite{JKSCounterBehrend}. Also, there is a more recent result, \cite{PSU26}, on the Behrend number of fat points in $\BA^N$ using Rees valuations.

\begin{remark}
    This work was developed completely independently of and concurrently with \cite{PSU26}, which appeared on arXiv in August 2026. While both studies independently derive the formulas for the zero-dimensional case, the present work specifically highlights the arbitrary-dimensional case. Also, while both studies reach similar conclusions for the zero-dimensional case, the different languages and tools used confirm the independence of the results.  After conversations with one of the authors, we noticed the following overlaps in the zero-dimensional case: Restricting the formulas for the general case and the monomial case (in Theorem  \ref{generalConjecture} and Theorem \ref{gcdMainThm}, respectively) to dimension zero are equivalent to the formulas in \cite{PSU26}. Moreover,  the formula for the monomial case in terms of the constant term of the defining hyperplanes in Theorem \ref{ThmC}, restricted to dimension zero, is equivalent to the results in \cite{PSU26}. 

\end{remark}

 \subsection*{Convention}  
 
By a \textit{polytope}, we mean a convex bounded polytope. Furthermore, we always consider ideals with a minimal number of generators. By a \textit{lattice point}, we mean a point with integer coordinates. Also, a \textit{lattice polytope} is a polytope such that all its vertices are lattice points. For simplicity, in the polynomial ring with two variables we use $x,y$ instead of $x_1,x_2$, and in the polynomial ring with three variables we use $x,y,z$ instead of $x_1,x_2,x_3$. A generator of a monomial ideal in $\BC[x_1,\ldots,x_N]$ is a \textit{mixed term} or a \textit{mixed generator}, if it has more than one variable $x_i$. For example $x^2y$ is a mixed generator of the ideal $(x^3,z^3,x^2y)$ in $\BC[x,y,z]$. We use \emph{the convex hull of an ideal} and \emph{the convex hull of the normal cone of the ideal} interchangeably. We denote the origin of $\BC^N$ by $\mathbf{o}$. The proofs are considered for $\mathbf{o}$, but it can be replaced by any singular point of $X$. When $X$ is zero-dimensional, we simply denote $v_X(\mathbf{o})$ by $v_X$. As mentioned above, the term degree will be used for either ideal/quotient rings or their corresponding affine or (quasi-)projective schemes.

\subsection*{Notation} For an $r$-dimensional ideal $I$ in $A=\BC[x_1,\ldots,x_N]$ and $X=\Spec(A/I)$ with the origin $\mathbf{o}$ as its singular point, we have the following notations for some numbers/notions associated with $I$ (or $X$). For the definition of degree, see Definition \ref{def: degree}.

\begin{center}
   \begin{tabular}{ r l }
     $\CN(I)$ & The normal cone of $I$ (see Definition \ref{Def:normalCone}).\\
        
      $\CS(I)$ & The special fibre of $I$ as the reduction of $\CN(I)$
      (see Definition \ref{Def:normalCone}).
      \\

      $n(I)=\sum_{j=0}^{r}n_j(I)$ & Number of the maximal dimensional irreducible components\\& of  $\CN(I)$ or $\CS(I)$ (see definition \ref{Def:normalCone}), where $n_j(I)$ is the\\& number of the irreducible components of $\CN(I)$ of the maximum\\& dimension, dominating a $j$-dimensional component of $\Spec(A/I)$.\\
       $\CN_i$ or $\CN_{i,j}$ & Components of  ${\CN(I)}$, for $I$ zero-dimensional or higher dimensional,\\& respectively (see Definition \ref{Def:normalCone}).\\
         $\mult(I)$ &The normalised leading coefficient of the Hilbert polynomial of $A/I$.\\ 

       $\mult_{i,j}$&  The degree of the $p_{i,j}$-primary decomposition component $Q_{i,j}$ of the\\& ideal of the normal cone $\CN(I)$. When $I$ is zero-dimensional, we\\& simply denote this by $\mult_i$.\\

       $D_{i,j}$& The degree of $p_{i,j}=\sqrt{Q_{i,j}}$, where $Q_{i,j}$ are primary decomposition\\& components of the ideal of the normal cone $\CN(I)$ for $0\leq j \leq r$, and \\& $1\leq i \leq n_j(I)$. When $I$ iszero-dimensional, we simply denote this by $D_i$.\\

          $l(\CN_{i,j})$& The length of $\CN(I)$ at the generic point of $\CN_{i,j}$ (see Definition \ref{def: distinguished}).\\

  $\mathfrak{c}_{X}$&Integral cycle on $X$ (see Definition \ref{def: distinguished}).\\
$\mathrm{Eu}$&Local Euler obstruction (see \cref{Def?:Behrend2}).\\
       
  \end{tabular}
   \end{center}
 \begin{center}
 \begin{tabular}{r l}

    $v_X(\mathbf{o})$ & Behrend function of $X$ at the origin, when $X$ is positive dimensional. \\&(see Definition \ref{Def: Behrend} and Definition \ref{Def?:Behrend2}).\\
    $v_X$ & Behrend function (Behrend number) of $X$, when $X$ is zero-dimensional.\\
    &(see Definition \ref{Def: Behrend} and Definition \ref{Def?:Behrend2}).\\
     
   \end{tabular}
\end{center}

\subsection*{Acknowledgements} 

The author thanks Michele Graffeo, Joachim Jelisiejew, and Andrea Ricolfi, for several discussions at the beginning of the project in 2021/2022. The author also thanks Hal Schenck for several helpful discussions in 2024.  The author is grateful to Mark Gross for his support over the years. The author benefited from discussions with David Eisenbud, Barbara Fantechi, Rahul Pandharipande, 
 Francisco Santos,  Johannes Schmitt, Mike Stillman, Richard Thomas, Ajith Urundolil-Kumaran,
Ravi Vakil, Andrzej Weber
 and Volkmar Welker. We thank Yunfeng Jiang, Andrea Ricolfi, Alessio Sammartano and Richard Thomas for comments on the final draft. This material is partially based upon work supported by the NSF under Grant No. DMS-1928930, while the author was in residence at SLMath (formerly known as MSRI) in Berkeley during the Fall
2022. The author was supported by the ERC Advanced Grant MSAG, the UKRI grant No. EP/X032779/1, and was partly in residence at ETH Z\"urich. The author thanks the Isaac Newton Institute for providing the space to discuss with the members of the program ``New Equivariant Methods in Algebraic and Differential Geometry''.  The rough ideas of the general subject were initiated when the author was partially supported by EPSRC grants EP/T015896/1 and EP/W522648/1.
The author also acknowledges using GeoGebra (\cite{geogebra5}) and Macaulay2 (\cite{GS}). The codes in Appendix \ref{App: M2Codes} are duo to Mike Stillman, while visiting Cambridge in 2023. While searching for relevant recent results, we noticed a recent preprint, \cite{PSU26}, on the zero-dimensional case; we acknowledge that our work is independent of it.

\section{Definitions and Commutative algebra background}
 In this section, we give some basic definitions related to the original definition of the Behrend function. We also review some definitions and facts in commutative algebra that we need in the following sections. We first define the Behrend function as in \cite[Definition 3.6, and Introduction]{Beh}.

\begin{definition}[Degree of an ideal or a qoutient ring]\label{def: degree}
    Throughout the paper, by \textit{degree} of an ideal or a quotient ring, we mean the factorial of the (Krull) dimension of the ideal or the  quotient ring, times the leading coefficient of the corresponding Hilbert polynomial.

    We also use the term `degree' for the corresponding affine or (quasi-)projective scheme.
\end{definition}

 \subsection{Behrend function}\label{subsection:Behrend} We have the following original key definition.fire
 \begin{definition}[Behrend function, \cite{Beh}]\label{Def: Behrend}
    Let $X$ be a scheme over complex numbers. The \textit{Behrend function} on $X$ is defined to be a constructible function 
    \begin{align*}
        v_X:X\to\BZ,
    \end{align*}
    so that if the scheme is proper and embeddable, then for  any symmetric obstruction theory on it with associated  \textit{virtual count} $\deg([X]^{\mathrm{vir}}):=\int_{[X]^{\mathrm{vir}}}1$, we have
\begin{align}\label{eq: WeightedEuler}
    \int_{[X]^{\mathrm{vir}}}1=\chi(X, v_X),
\end{align}
where $[X]^{\mathrm{vir}}$ is the virtual fundamental class, and $\chi(X, v_X)$ is the weighted Euler characteristic.
 \end{definition}

To give a more precise definition of the Behrend function, we need to define the zero locus cycle associated with $X$. We use the version in \cite[Section 1.1]{Beh}.
 \begin{definition}[Distinguished cycle, \cite{Beh}]\label{def: distinguished}
     For an ideal $I$ in $A=\BC[x_1,\ldots,x_N]$, let $\CN(I)$ be the normal cone of $I$, with irreducible components of maximal dimension $\CN_{i,j}$, and let $X=\Spec(A/I)$. Let $l(\CN_{i,j})$ be the length of the normal cone at the generic point of $\CN_{i,j}$.

     Then, there is a canonical projection $\pi: \CN(I) \to X$. Now define 
     \begin{align}\label{eq:distinguished}
        \mathfrak{c}_{X}=\mathfrak{c}_{X/{\BA^N}}:=\sum_{\CN_{i,j}\subset\CN(I)}(-1)^{\dim\pi(\CN_{i,j})}l(\CN_{i,j})\big[\pi(\CN_{i,j})\big],
     \end{align}
     to be the unique integral \textit{cycle on $X$} (or the \emph{distinguished cycle on $X$}). 
 \end{definition}
We also need the notion of local Euler obstruction as originally defined in \cite[Section 3]{MacPherson}. We use the version defined in \cite[Section 1.2]{Beh}.
 \begin{definition}[Local Euler obstruction, \cite{Beh,Jiang}]\label{Def?:Behrend2} Let $X$ be as in Definition \ref{def: distinguished}. We define
 \begin{align*}
     \mathrm{Eu} : Z_*(X)&\to \mathrm{Cons}(X)\\
     \{\text{integral algebraic cycles on $X$}\} &\mapsto \{\text{constructible $\BZ$-valued
functions on $X$}\}.
 \end{align*}
 More precisely, let $X\hookrightarrow \mathbb{A}^N$ be an embedding and $\dim X=r$. Let $(x_1,\ldots,x_N)$ be local coordinate at $p\in X$. Also, define $X^{\circ}$ to be the smooth locus in $X$ and $s:X^{\circ}\to \mathrm{Gr}_r(T\mathbb{A^N})$ be a section defined by $s(x)=T_x(X^{\circ})$, where $\mathrm{Gr}_r(T\mathbb{A^N})$ is the Grassmannian bundle of $r$ dimensional subspaces of $T\mathbb{A}^N$ over $\mathbb{A}^N$. Then, we can define $\hat{X}$ to be the closure of the image of $s$. Let $\alpha:\hat{X} \to X$, defined by restriction. Also, define the function $\norm{X}$ to be $\sqrt{x_1\bar{x_1}+\ldots+x_N\bar{x_N}}$; so, the differential $\mathrm{d}\norm{X}$  is a section of $T^{*}\mathbb{A}^N$. Let $t$ be the pullback of  $\mathrm{d}\norm{X}$ under the embedding. Define $B_{\varepsilon}:=\{x:\norm{x}\leq \varepsilon\}$ and $S_{\varepsilon}:=\{x:\norm{x}=\varepsilon\}$. We can then define
 \begin{align*}
     \mathrm{Eu}(X)(p):=\langle\mathrm{Eu}(T^*X,t),\mathcal{O}\rangle,
 \end{align*}
     where $\mathcal{O}\in H_{2r}(\alpha^{-1}B_{\varepsilon},\alpha^{-1}S_{\varepsilon};\mathbb{Z} )$ and $\mathrm{Eu}(T^*X,t)$ denotes the
obstruction to extending $t$ as a non-zero section of $T^*X$ from $\alpha^{-1}S_{\varepsilon}$ to $\alpha^{-1}B_{\varepsilon}$.
 \end{definition}
Now, the Behrend function of $X$ as above (Definition \ref{def: distinguished}) can be defined as (See \cite[Definition 1.4]{Beh})
\begin{align}\label{eq:BehEuler}
    v_X:=\mathrm{Eu}(\mathfrak{c}_X).
\end{align}

Applying $\mathrm{Eu}$ to \ref{eq:distinguished}, we get 
 \begin{align}
       v_X=\sum_{\CN_{i,j}\subset\CN(I)}(-1)^{\dim\pi(\CN_{i,j})}l(\CN_{i,j})\mathrm{Eu}(\big[\pi(\CN_{i,j})\big]),
     \end{align}

     We have the following properties for.
     \begin{lemma}\label{PropertiesEuler}\cite[Definition 2.4]{Jiang},\cite[Section 3]{MacPherson}, \cite{RicolfiSign} We have
     \begin{itemize}
         \item  $\mathrm{Eu}(X)(p) = 1$ if $Z$ is nonsingular at $p$.
         \item $\mathrm{Eu}(X)(p) = 0$ if $p\not\in Z$
         \item$ \mathrm{Eu}(X_1 \times X_2)(p_1 \times p_2) = \mathrm{Eu}(X_1)(p_1) \cdot \mathrm{Eu}(X_2)(p_2)$
\item If  $X$ or Deligne-Mumford stack which is reducible at the point $p$, and
let $X_i$ are the irreducible components, then $\mathrm{Eu}(X)(p) =\sum\mathrm{Eu}(X_i)(p)$.
     \end{itemize}

     \end{lemma}

We give another definition for future reference.
\begin{definition}[Dominating the support]\label{def: dominant} We say that a maximum-dimensional irreducible component $\CN_{i,j}$ of $\CN(I)$ \emph{dominates}  the support $\pi(\CN_{i,j})$ in $X$.
    
\end{definition}
\begin{remark}
    Note that on a smooth point of $X$, the value of $v_X$ equals $(-1)^{\dim X}$. Thus, the value of the Behrend function is interesting only on the singular points of $X$, which is the subject of this article.
\end{remark}

If $X$ is a critical locus, its Behrend function has a particular enumerative geometric interpretation as in \eqref{eq: WeightedEuler}, and, in theory, can be computed as follows.
\begin{lemma}\cite[Section 1.2]{Beh}\label{lem: CriticalLocus}
If $X$ is a critical locus, i.e., if $X$ is the zero locus of $df$, for some regular function $f$ on $\BA^N$, then for any $\BC$-valued point $p\in X$, we have
\begin{align}
    v_X(p)=(-1)^{N}(1-\chi(F_p)),
\end{align}
where $F_p$ the Milnor fibre of $f$ at $p$.
\end{lemma}
\begin{proof}
    See \cite[Corollary 2.4 (iii)]{PP}.
\end{proof}

 We devote the rest of the section to the commutative algebraic definitions.
\subsection{Associated primes and primary decomposition}

\begin{definition}\cite[Chapter 4]{Atiyah-MacDonaldMacDonald} Let $R$ be a ring. An ideal $q$ in $R$ is \textit{primary} if it is not equal to the whole ring and if $xy\in q$, then either $x\in q$ or $y^m\in q$, for some $m\in\BZ^{\geq 1}$. If $p=\sqrt q$, then $q$ is called \textit{$p$-primary}.
\end{definition}
\begin{definition}\cite[Section 3.1]{Eisenbudcomm} Let $R$ be a ring and $M$ be an $R$-module. A prime ideal $p$ in $R$ is called \textit{associated} to $M$ if it is the annihilator of some $x\in M$.    
\end{definition}
The set of all prime ideals associated with $M$ is denoted by $\mathrm{Ass}_R(M)$ or $\mathrm{Ass}(M)$.
\begin{lemma}\cite[Theorem 3.1]{Eisenbudcomm} Let $R$ be a Noetherian ring and $M\neq0$ be a finitely generated $R$-module. Then, $\mathrm{Ass}(M)$ is a finite non-empty set of primes, each containing the annihilator of $M$. The set $\mathrm{Ass}(M)$ includes all primes minimal among primes containing $\mathrm{ann}(M)$. 
    
\end{lemma}

If $R$ is a commutative ring, the minimal elements in $\mathrm{Ass}_R(M)$ are called \textit{isolated primes}, and the non-minimal ones are called \textit{embedded primes}.

\subsection{Noether normalization}
The following lemma, which is called \textit{Noether normalization lemma}, is a version of it that we need later:
\begin{lemma}\cite[Lemma 10.115.4]{stacks-project} \label{noetherNormalization}
    For a field $k$, and an $m$-dimensional ideal $I$, let $S=\frac{k[x_1,\ldots,x_N]}{I}$. Then, there exist a polynomial ring $k[y_1,\ldots,y_m]$ and a finite injective homomorphism $k[y_1,\ldots,y_m] \hookrightarrow S$.
\end{lemma}

\subsection{Normalization and convex hull}

\begin{definition}
Let $I$ be a monomial ideal in $\BC[x_1,\ldots,x_N]$. An element $b\in \BC[x_1,\ldots,x_N]$ is \textit{integral} if there are $c_i\in I^{i}$ such that 
\begin{align*}
    b^{l}+c_1b^{l-1}+\ldots+c_l=0.
\end{align*}

    The \textit{integral closure} of a monomial ideal $I$ is denoted by $\bar{I}$ and consists of all the integral elements as above. A monomial ideal $I$ is  \textit{integrally closed} if $\bar{I}=I$.
\end{definition}

\begin{definition}
    We define the geometric \textit{normal cone} of an ideal $I$ in $\BC[x_1,\ldots,x_N]$ to be $\Spec(\mathrm{gr}(I))=\Spec(\bigoplus_{i=0}^{\infty}\frac{I^i}{I^{i+1}})$. Note that we may use the word normal cone for either $\Spec(\bigoplus_{i=0}^{\infty}\frac{I^i}{I^{i+1}})$ or $\bigoplus_{i=0}^{\infty}\frac{I^i}{I^{i+1}}$, interchangeably.
\end{definition}
Note that minimal primes of  $\mathscr{I}({\CN(I)})$, the defining ideal of $\CN(I)$, are equivalent to the maximal dimensional components of the normal cone $\CN(I)$.

\begin{lemma}\cite[Exercise 4.23]{Eisenbudcomm}, \cite[Proposition 1.4.6]{Huneke-Swanson} \label{normalizationConvexHull} Let $I$ be a monomial ideal in $\BC[x_1,\ldots,x_N]$, and $\text{\:}\Gamma$ be the set of exponents of monomials in $I$. Let $\Lambda$ be the convex hull of $\text{\:}\Gamma$, and $\Lambda^*$ be the set of integral points in  $\Lambda$. Then, the integral closure of $I$ is the ideal generated by $\Lambda^*$.
\end{lemma}

\begin{corollary} \label{cor: restriction}
    Let $I$  be a monomial integrally closed ideal in $\BC[x_1,\ldots,x_N]$. Then, the restriction of $I$ to any non-empty subset of variables $\{x_1,\ldots,x_N\}$ is still integrally closed.
\end{corollary}

\begin{proof}
 By Lemma \ref{normalizationConvexHull}, there is a one-to-one correspondence between the normal cone of an ideal and the convex hull of the exponents. As the restriction of the convex hull to any non-empty subset of variables $\{x_1,\ldots,x_N\}$ is still a convex hull, the same lemma implies the result.
\end{proof}

\section{The quotient formula: proof of the abstract version}

In this section, we give a simple proof of Theorem \ref{generalConjecture} (see Theorem  \ref{v as quotien a of mult (reducedness condition)} below) using algebraic arguments.
Then, in the following sections, we will see a geometric interpretation of this abstract version of the quotient formula.

\begin{theorem}
 \label{v as quotien a of mult (reducedness condition)}
Let $I$ be an 
ideal in $A=\BC[x_1,...,x_N]$ with $l$ minimum number of generators, and let $X=\Spec(A/I)$. Then,  for any $r$-dimensional ideal $I$ in $\BC[x_1,\ldots,x_N]$, and $\mathbf{o}$ the origin as a singular point of  $X=\Spec(A/I)$, we have
\begin{align*}
v_X(\mathbf{o})=&\sum_{j=0}^{r}\sum_{i=1}^{n_j(I)}(-1)^j\frac{\mult_{i,j}}{D_{i,j}}\mathrm{Eu}(\big[\pi(\CN_{i,j})\big]),
\end{align*}
where $\mult_{i,j}$ is the degree of the $p_{i,j}$-primary decomposition component of maximal dimension, $Q_{i,j}$, of $\CN(I)$, and $D_{i,j}$ is the degree of $p_{i,j}$, for $0\leq j \leq r$, and $1\leq i \leq n_j(I)$. Also, $\CN_{i,j}$ are the maximal-dimensional irreducible components of $\CN(I)$, and $\pi:\CN(I)\to X$ is the canonical projection, and $\:{}\mathrm{Eu}$ is the local Euler obstruction. Also, $n(I)=\sum_{j=0}^{r}n_j(I)$ is the number of the irreducible components of the maximal dimension of $\CN(I)$, and $n_j(I)$ is the number of the irreducible components of the maximal dimension of $\CN(I)$, dominating an $j$-dimensional component of $\Spec(A/I)$ in $\BA^N$.

\end{theorem}

\begin{proof}

Let $R=\bigoplus_{i=0}^{\infty}\frac{I^i}{I^{i+1}}$, and ${\mathcal{N}(I)}=\Spec (R)$.
We can write $R$ as $\frac{\BC[x_1,\ldots,x_N,w_1,\ldots,w_l]}{L}$, where $L=I+J$, and $J$ is an ideal in $\BC[x_1,\ldots,x_N,w_1,\ldots,w_l]$, and $l$ is the (minimum) number of generators of $I$. We can write $L=\cap_{i,j} Q_{i,j}$, where $\{Q_{i,j}\}$ is its $p_{i,j}$-primary decomposition. 

Let $\{\mathcal{N}_{i,j}\}_{\substack{0\leq j \leq r\\1\leq i \leq n_j(I)}}$ be the set of maximum-dimensional irreducible components of ${\mathcal{N}(I)}$, corresponding to $Q_{i,j}$. Think of $\CS(I)$ as the reduction of $\CN(I)$, and let  $\{\mathcal{S}_{\mathcal{N}_{i,j}}\}_{\substack{0\leq j \leq r\\1\leq i \leq n_j(I)}}$ be the set of corresponding irreducible components of  $\mathcal{S}(I)$. Since $\sqrt{Q_{i,j}}=p_{i,j}$, the components $\mathcal{S}_{\mathcal{N}_{i,j}}$ correspond to $p_{i,j}$.
Since $\dim(L)=N$, from Noether normalization lemma (Lemma \ref{noetherNormalization}), there is a finite injective map 
\begin{align}\label{eq:iota}
    \iota:\BC[y_1,\ldots,y_N]\hookrightarrow R.
\end{align}

For $F=\BC[y_1,\ldots,y_N]^{*}$, $\iota$ induces a free injection 
\begin{align}\label{eq:FinverseIota}
    F^{-1}\iota: \BC(y_1,\ldots,y_N) \hookrightarrow F^{-1}R,
\end{align}
where $\BC(y_1,\ldots,y_N)$ is the function field of $\BC[y_1,\ldots,y_N]$. Taking $\Spec$ of the finite map $\iota$ reverses it to 

$$\phi:\CN(I)\to \BA^N,$$
and taking $\Spec$ of $F^{-1}\iota$ gives 

$$\Spec(F^{-1}R)\to\Spec(\BC(y_1,\ldots,y_N)).$$ 

 Hence, the generic fibre of $\phi$  is given by $\Spec(F^{-1}R)$, which has irreducible components corresponding to maximum-dimensional ($=N=\dim(R)$) irreducible components of $\CN(I)$, i.e., $\CN_{i,j}$, for $0\leq j \leq r$ and $1\leq i \leq n_j(I)$.  
Therefore, the rank of the free map $F^{-1}\iota$ is $\mult(I)=\sum_{\substack{0\leq j \leq r\\1\leq i \leq n_j(I)}}\mult_{i,j}$.

Now, let $R_{i,j}=R_{p_{i,j}}=\CO_{\CN_{i,j},\CN(I)}$ be the local ring of $R$ at the defining ideal of $\mathcal{N}_{i,j}$; hence, $F^{-1}R$ can be written as $\prod_{\substack{0\leq j \leq r\\1\leq i \leq n_j(I)}}R_{i,j}$. 
Therefore, by taking the degree of the map $F^{-1}\iota$, we can write the multiplicity, $\mult(I)$, as $\sum_{\substack{0\leq j \leq r\\1\leq i \leq n_j(I)}}\dim_{\BC(y_1,\ldots,y_N)}R_{i,j}$.

On the other hand, 
 $\mathcal{S}(I)$ can be written as $\prod_{\substack{0\leq j \leq r\\1\leq i \leq n_j(I)}}K(R_{i,j})$, where $K(R_{i,j})$ is the residue field of $R_{i,j}$ and is a finite extension over $\BC(y_1,\ldots,y_N)$. 

Since $\BC(y_1,\ldots,y_N)$ is not necessarily algebraically closed, it is not necessarily equal to its finite extension.  The difference measures the difference between the multiplicity and the Behrend function as below.

Similarly as the case of multiplicity, the degree of $\CS(I)$ can be written as $$\sum_{\substack{0\leq j \leq r\\1\leq i \leq n_j(I)}}\dim_{\BC(y_1,\ldots,y_N)}K(R_{i,j}).$$ Note that by taking the residue field of $\iota$,
we obtain \begin{align} \label{eq: degreeResidue}
\deg(p_{i,j})
=\dim_{\BC(y_1,\ldots,y_N)}K(R_{i,j}).
\end{align}

Note that we can write $Q_{i,j}R_{i,j}=p_{i,j}^{n_{i,j}}R_{i,j}$, for some integers $n_{i,j}$. From this, we have $\deg(Q_{i,j})=\deg(p_{i,j})\length(\frac{R_{i,j}}{Q_{i,j}R_{i,j}})$. Thus, $\deg(p_{i,j})|\deg(Q_{i,j})$.

Recall that, by our definition, $\mult_{i,j}$ is the degree of  ${\CN_{i,j}}$, and $D_{i,j}$ is the degree  of $\CS_{\CN_{i,j}}$; i.e., $\mult_{i,j}$ is the contribution of $\mult(I)$ to ${\CN_{i,j}}$, and $D_{i,j}$ is the contribution of the degree of $\CS(I)$ to $\CS_{\CN_{i,j}}$.
 Now, from \eqref{eq: degreeResidue}, for each fixed $i,j$, for the component $\CS_{\CN_{i,j}}\subset \CN_{i,j}$, we have
\begin{equation}
    \begin{aligned} \label{eqn: multD}
\frac{\mult_{i,j}}{D_{i,j}}=&\frac{\mult|_{\CN_{i,j}}}{D|_{\CS_{\CN_{i,j}}}}=\frac{\deg{\CN_{i,j}}}{\deg{\CS_{\CN_{i,j}}}}=\frac{\deg(Q_{i,j})}{\deg(p_{i,j})}
\\
=&\frac{\dim_{\BC(y_1,\ldots,y_N)}R_{i,j}}{\dim_{\BC(y_1,\ldots,y_N)}K(R_{i,j})}=\dim_{K(R_{i,j})}R_{i,j}=\deg \CO_{\CN_{i,j},\CN(I)}=l(\CN_{i,j}),
\end{aligned}
\end{equation}
where $l(\CN_{i,j})$ is the length of $\CN(I)$ at (the generic point of) $\CN_{i,j}$. 

Restricting to $r=0$, we have $n_0(I)=n(I)$. Also, in this case, we have the notation: $\mult_i:=\mult_{i,0}$, $D_i:=D_{i,0}$ $R_i:=R_{i,0}$, $Q_i:=Q_{i,0}$, $p_i:=p_{i,0}$, $\CN_i:=\CN_{i,0}$, and  $\CS_{\CN_i}:=\CS_{\CN_{i,0}}$. Then, \eqref{eqn: multD} becomes
\begin{align} \label{eqn: multD0}
\mathrm{dim}_{K(R_i)}R_i=\deg \CO_{\CN_{i},\CN(I)}=l(\CN_{i})=v_{X,i},
\end{align}
where $v_{X,i}$ is the contribution of $v_X$ to each irreducible component of the normal cone, and the last equality comes from the definition of the Behrend number for the zero-dimensional case (see in  \eqref{eq:distinguished} and \eqref{eq:BehEuler}).  Hence, for zero-dimensional $I$ (for which, the contribution of $\mathrm{Eu}$ is $1$ by Lemma \ref{PropertiesEuler}), from \eqref{eqn: multD0}, we have
\begin{align}\label{eqn:0dim}
\sum_{i=1}^{n(I)}\frac{\mult_i}{D_i}
=\sum_{i=1}^{n(I)}v_i(I)=v_X.
\end{align}

Now, let $r>0$, and $X=\Spec(A/I)$ with a singular point at the origin $\mathbf{o}\in X$.
Following Definition \ref{def: distinguished}, let us consider the canonical projection $\pi: \CN(I) \to X$. By definition \ref{def: dominant}, any maximum-dimensional irreducible component of the normal cone, $\CN_{i,j}\subset \CN(I)$, dominates its corresponding $i$-dimensional support $\pi(\CN_{i,j})\subset X$, for each $0\leq i \leq r$ and fixed $1\leq j \leq n_i(I)$. Note that since we assume that the origin is a singular point for $X$, the components of $X$ are coordinate planes of possibly different dimensions (as in Remark \ref{coordinatePlanes}).  Hence, the general alternating sum in the formula comes from the formula for the cycle on $X$ introduced in Definition \ref{def: distinguished} (also, see \cite[Section 2.2]{RicolfiSign}), and also  \eqref{eqn: multD} and \eqref{eqn:0dim} above. This means that we have
\begin{align*}
\sum_{j=0}^{r} \sum_{i=1}^{n_j(I)}(-1)^j\frac{\mult_{i,j}}{D_{i,j}}\mathrm{Eu}(\big[\pi(\CN_{i,j})\big])=&\sum_{j=0}^{r}\sum_{i=1}^{n_j(I)}(-1)^j l(\CN_{i,j})\mathrm{Eu}(\big[\pi(\CN_{i,j})\big])=  v_X(\mathbf{o}),
\end{align*}
as claimed.

\end{proof}

\begin{corollary}[Upper bound on the Behrend number]
For any zero-dimensional ideal $I$ in $A=\BC[x_1,...,x_N]$ and $X=\Spec(A/I)$, we have $$v_X\leq \mult (I).$$
\end{corollary}
\begin{proof}
Note that $\mult(I)=\sum_{i=1}^{n(I)}\mult_i(I)$. So by Theorem \ref{v as quotien a of mult (reducedness condition)}(iii), we obtain the result.
\end{proof}

\section{Minimisation}\label{minimisation} 
In this section, we introduce the notion of minimisation and apply it to give a more practical interpretation of the quotient formula for monomial ideals. 
This will also characterise monomial ideals for which the Behrend number is maximal (see Proposition \ref{maxB}). 

\begin{convention*}
    In this section, we assume that
the ideals have a minimum number of generators. 
\end{convention*}

\begin{definition}[Minimal ideal]
We say that a monomial ideal $I$ with a minimum number of generators in $\BC[x_1,x_2,\ldots,x_N]$ is \textit{minimal} if the removal of any mixed generator from it changes the multiplicity, that is, if $I= (x_1^{m_1}, x_2^{m_2},\ldots,x_N^{m_N},F_1,F_2,\ldots , F_k)$, where $m_i$ could possibly be zero (if $I$ is of a dimension greater than zero),and $F_i$ are mixed monomial terms. Then, $$\mult (x_1^{m_1}, x_2^{m_2},\ldots,x_N^{m_N},F_1,F_2,\ldots , F_k)\neq \mult (x_1^{m_1}, x_2^{m_2},\ldots,x_N^{m_N},F_1,F_2,\ldots,\widehat{F_i},\ldots, F_k),$$ for any omitted mixed term $\widehat{F_i}$. The \textit{minimisation} of $I$, which is indicated by $I_{\min}$, is a minimal ideal with a minimal number of generators such that the set of its generators is contained in the set of the (minimal) generators of $I$ such that $\mult (I)=\mult (I_{\min})$.

Equivalently, an ideal $I$ is \textit{minimal}, if and only if $I=I_{\min}$.
\end{definition}

The following lemma is a consequence of the definition:

\begin{lemma}
Let $I=(I_{\min}, H_1,H_2,\ldots,H_l)$. Then
$$\mult (I)=\mult (I_{\min}, H_1,H_2,\ldots,H_{l-1})=\mult (I_{\min}, H_1,H_2,\ldots,H_{l-2})=\ldots=\mult (I_{\min}).$$
\end{lemma}

Now, we list some properties of minimisation in the following proposition, which are direct consequences of the combinatorial interpretation of minimisation in the following sections:

\begin{prop}\label{nImin}Let  $I$ be a monomial ideal in $\BC[x_1,x_2,\ldots,x_N]$.\begin{itemize}
\item  Minimisation preserves normalization: 
$ \text{normalization of $I$= normalization of $I_{\min}$}$.

    \item  Normalization preserves minimisation: $\text{$I_{\min}=($normalization of $I)_{\min}$}$.

\item Minimisation preserves the number of maximal dimensional irreducible components of the normal cone of the ideal: $n(I)=n(I_{\min})$.
\end{itemize}
   
\end{prop}

\subsection{An algorithm for minimisation} \label{minalgorithm}

Let $I=(x_{1}^{m_1},x_{2}^{m_2},\ldots,x_{N}^{m_N},N_1,N_2,\ldots,N_s,F_1,F_2,\ldots F_{h})$ be an ideal minimally generated by monomials in $\BC[x_1,x_2,\ldots,x_N]$, such that $F_i,N_i$ are mixed terms and $m_i$ possibly being zero for some (or all) $i$. Also, suppose that $N_i$ is contained in the (minimal set of generators of the) normalization of $I$ for $1\leq i\leq s$, and $F_j$ is not in the (minimal set of generators of the) normalization of $I$ for $1\leq j\leq h$. Then, we have the following algorithm to find $I_{\min}$:

\begin{enumerate}
    \item Remove all $F_j$'s for $1\leq j\leq h$, so that we are left with $(x_{1}^{m_1},x_{2}^{m_2},\ldots,x_{N}^{m_N},N_1,N_2,\ldots,N_s)$,
    \item Omit each $N_i$ one by one from $I$ and check if $N_i$ appears in the (minimal set of generators of the) normalization of  $(x_{1}^{m_1},x_{2}^{m_2},\ldots,x_{N}^{m_N},N_1,N_2,\ldots,\Hat{N_i}, \ldots N_s)$.
    \begin{enumerate}
        \item if $N_j$ appears in the (minimal set of generators of the) normalization of  $$(x_{1}^{m_1},x_{2}^{m_2},\ldots,x_{N}^{m_N},N_1,N_2,\ldots,\Hat{N_j}, \ldots,N_s),$$ then mark $N_j$ to remove it after completing step (2) for each $N_i$.
        \item if $N_j$ does not appear in the (minimal set of generators of the) normalization of  $(x_{1}^{m_1},x_{2}^{m_2},\ldots,x_{N}^{m_N},N_1,N_2,\ldots,\Hat{N_j}, \ldots,N_s)$, then mark $N_j$ to keep it after completing step (2) for each $N_i$. We rename such $N_j$ to $N'_j$.
    \end{enumerate}
    \item The ideal generated by $x_{1}^{m_1},x_{2}^{m_2},\ldots,x_{N}^{m_N}$ and the monomials in (b) is the (minimal set of generators of the) minimisation of $I$; i.e., $I_{\min}=(x_{1}^{m_1},x_{2}^{m_2},\ldots,x_{N}^{m_N},N'_1,N'_2,\ldots,N'_{t})$, where $t \leq s$.
  
\end{enumerate}

\begin{remark}
Note that in step (2), ordering does not matter. We remove one $N_i$ at a time, check the process for it, and then return it and repeat the process for another $N_i$.
\end{remark}

Let us look at an example and apply the algorithm above:
 
\begin{example}
Let $$I=(x^{200},y^{200}, x^2y^{20},x^6y^{19},x^{11}y^{17},x^{13}y^{12},x^{15}y^{11},x^{16}y^{7},x^{21}y^{5},x^{22}y^{4},x^{24}y^{3},x^{35}y^{2})$$  be an ideal in $\BC[x,y]$. One can easily check that $$(F_1,F_2,\ldots,F_5)=(x^6y^{19},x^{11}y^{17},x^{13}y^{12},x^{15}y^{11},x^{21}y^{5}),$$
and $$(N_1,N_2,\ldots,N_5)=(x^2y^{20},x^{16}y^{7},x^{22}y^{4},x^{24}y^{3},x^{35}y^{2}),$$$$(N'_1,N'_2,\ldots,N'_4)=(x^2y^{20},x^{16}y^{7},x^{24}y^{3},x^{35}y^{2}).$$ Therefore, we have $I_{\min}=(x^{200},y^{200},x^2y^{20},x^{16}y^{7},x^{24}y^{3},x^{35}y^{2})$.
\end{example}

  \section{Combinatorial interpretation: background}

This section provides a combinatorial background for giving a geometric interpretation of the quotient formula for monomial ideals. The first step is to correspond the monomial terms of a monomial ideal in $\BC[x_1,\ldots,x_N]$ to points in $\BZ^N$; i.e., we have
\begin{align*}
    \text{monomial term $x_1^{m_1}x_2^{m_2}\ldots x_N^{m_N}$}\longleftrightarrow\text{point $(m_1,m_2,\ldots,m_N)$ in $\BZ^N$=\text{lattice point}},
\end{align*}
where some $m_i$ can possibly be zero.

\subsection[$(N-1)$-boundary]{A combinatorial interpretation of the minimisation: determination of vertices in $(N-1)$-boundary} We make a combinatorial sense of the abstract notion of minimisation introduced previously. 

\begin{definition}[$(N-1)$-boundary]
For a minimal ideal with a minimal number of generators \begin{align*}
(x_1^{m_1},\ldots,x_N^{m_N},x_1^{m_{11}}x_2^{m_{21}}\ldots x_N^{m_{N1}},\ldots,x_1^{m_{1k}}x_2^{m_{2k}}\ldots x_N^{m_{Nk}})\end{align*} in $\BC[x_1,\ldots,x_N]$, where some $m_i$ can possibly be zero, the \textit{$(N-1)$-boundary}  is the lower boundary of the unique convex hull generated by $m_i$ or $m_{ij}$,  
for $1\leq i \leq N$ and $1\leq j \leq k$. 
\end{definition}

The existence and uniqueness of such a decomposition passing through any given finite set of points (after omitting unsuitable points) in $\BR^N$ is a well-known combinatorial fact; for example, see \cite[Corollary 1.7]{bruns2009polytopes}.  The following proposition summarises this for zero-dimensional monomial ideals:

\begin{prop} \label{prop: simplexDecomposition}
    Let $I$ be a zero-dimensional monomial ideal in $\BC[x_1,\ldots,x_N]$, and let $P_1=(p_{11},p_{12},\ldots p_{1N})$,$P_2=(p_{21},p_{22},\ldots p_{2N}),\ldots,P_k=(p_{k1},p_{k2},\ldots p_{kN})$ be the corresponding points (vectors) in $\BR^N$ associated with the exponents of the generators of $I$. Then, for some $j\leq N$, there exists a convex $(j-1)$-boundary in  $(\BR^{\geq 0})^{j}$ with the maximum number of $(j-1)$-polytopes (as lower facets of the corresponding convex hull) 
so that:
\begin{itemize}
\item the vertices of the $(j-1)$-polytopes cover all $P_i$'s

    \item the vectors corresponding to the vertices of each $(j-1)$-simplex generate an $j$ dimensional parallelepiped.
\end{itemize}

\end{prop}

The following corollary shows that minimisation provides a way to present exactly the vertices of the $(N-1)$-boundary:

\begin{corollary}
Minimisation removes the non-vertices monomials in the corresponding $(N-1)$-boundary; i.e., minimal ideals only consist of vertices.
\end{corollary}
\begin{proof} 

   Without loss of generality, we consider the zero-dimensional monomial ideals in $2$ variables. The proof for the general case is similar.
   
   Let $I=(A,B,C,D,E)\subset \BC[x,y]$, be a monomial ideal, where $B=x^{m_1}, A=y^{m_2}$, as in Figure \ref{Fig:three cases}. 

First, we note that if a monomial generator is not among the minimal generators of the normalization (which is the case for the mixed monomial term $E$ in Figure \ref{Fig:three cases}), then it will not appear on the boundary of the convex hull of the ideal; therefore, if the generators omitted in step (1) of the minimisation are not among vertices. For step (2), if we temporarily ignore a generator $g$ that is among the minimal generators of normalization but is not in the normalization of the rest of the generators, then that generator must be a vertex (this is the case for vertex $C$ in Figure \ref{Fig:three cases}(ii)).

Finally, if the monomial is not a vertex but is on the lower boundary of the convex hull, then by definition it will be omitted in the minimisation process (This is the case for the mixed monomial generator $D$ in Figure \ref{Fig:three cases}(iii)).

\begin{figure}[h]
 \subcaptionbox*{(i)}[.30\linewidth]{%
    \includegraphics[width=\linewidth]{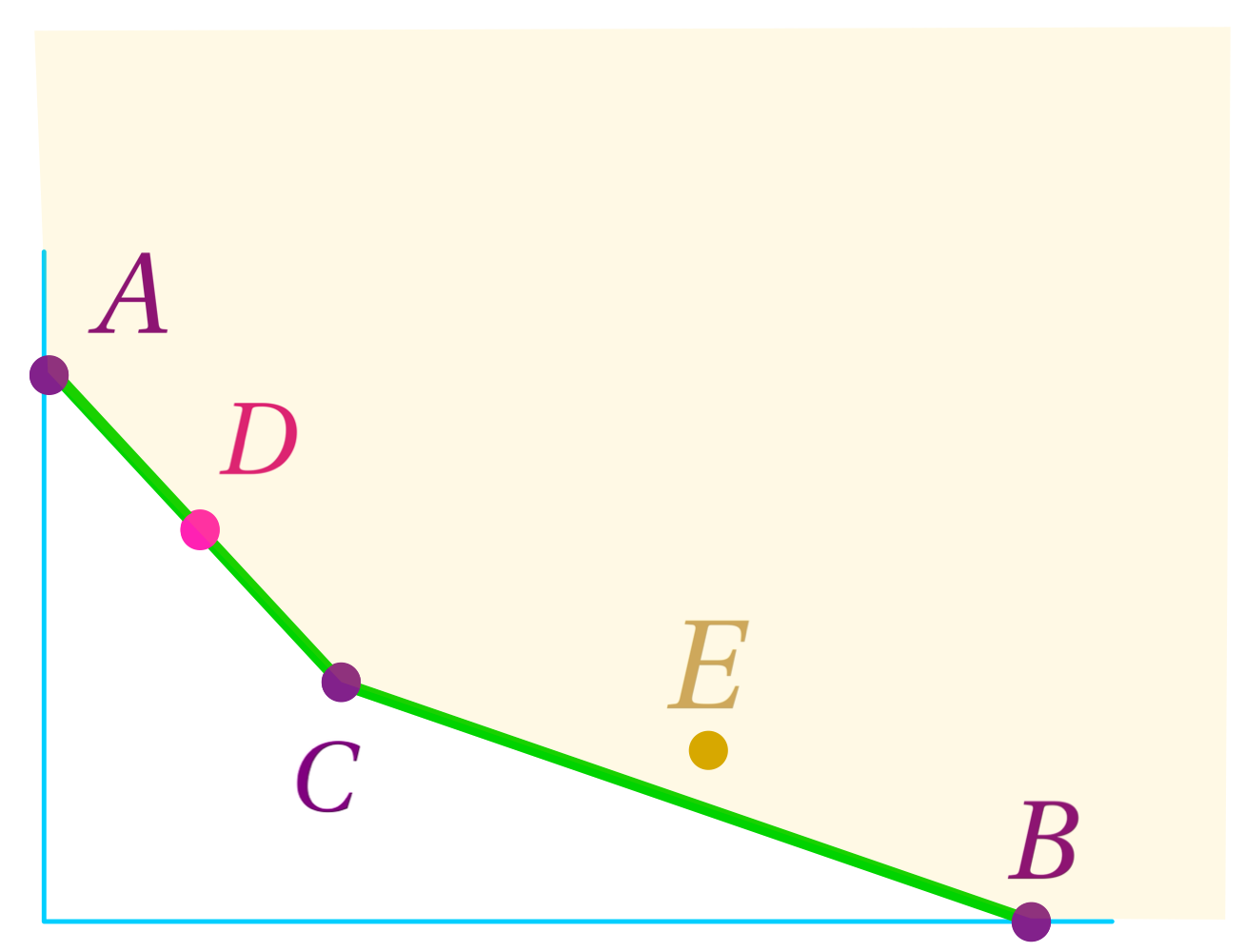}%
  }%
 \hskip3ex
  \subcaptionbox*{(ii)}[.30\linewidth]{%
    \includegraphics[width=\linewidth]{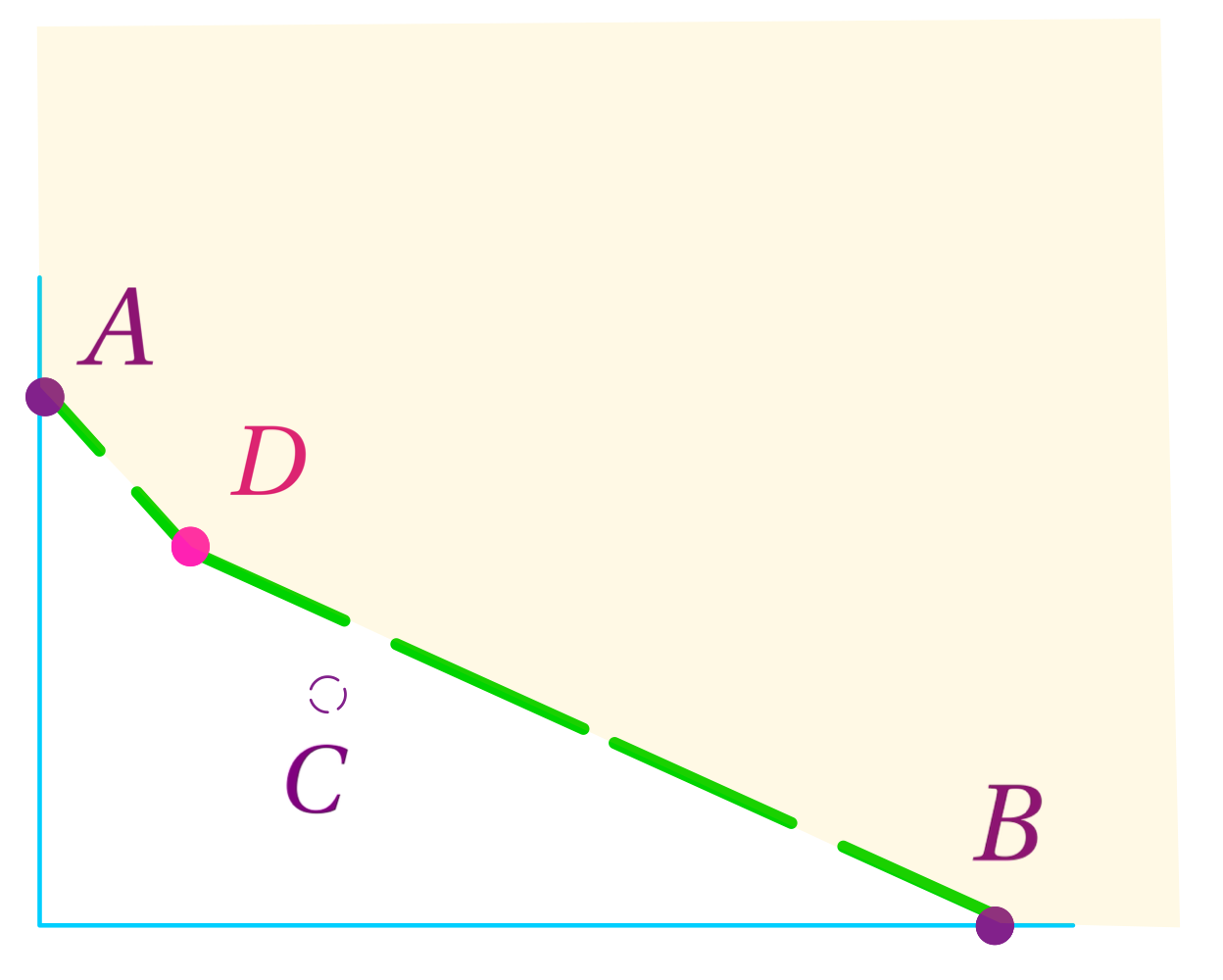}%
 }
  \hskip3ex
  \subcaptionbox*{(iii)}[.30\linewidth]{%
   \includegraphics[width=\linewidth]{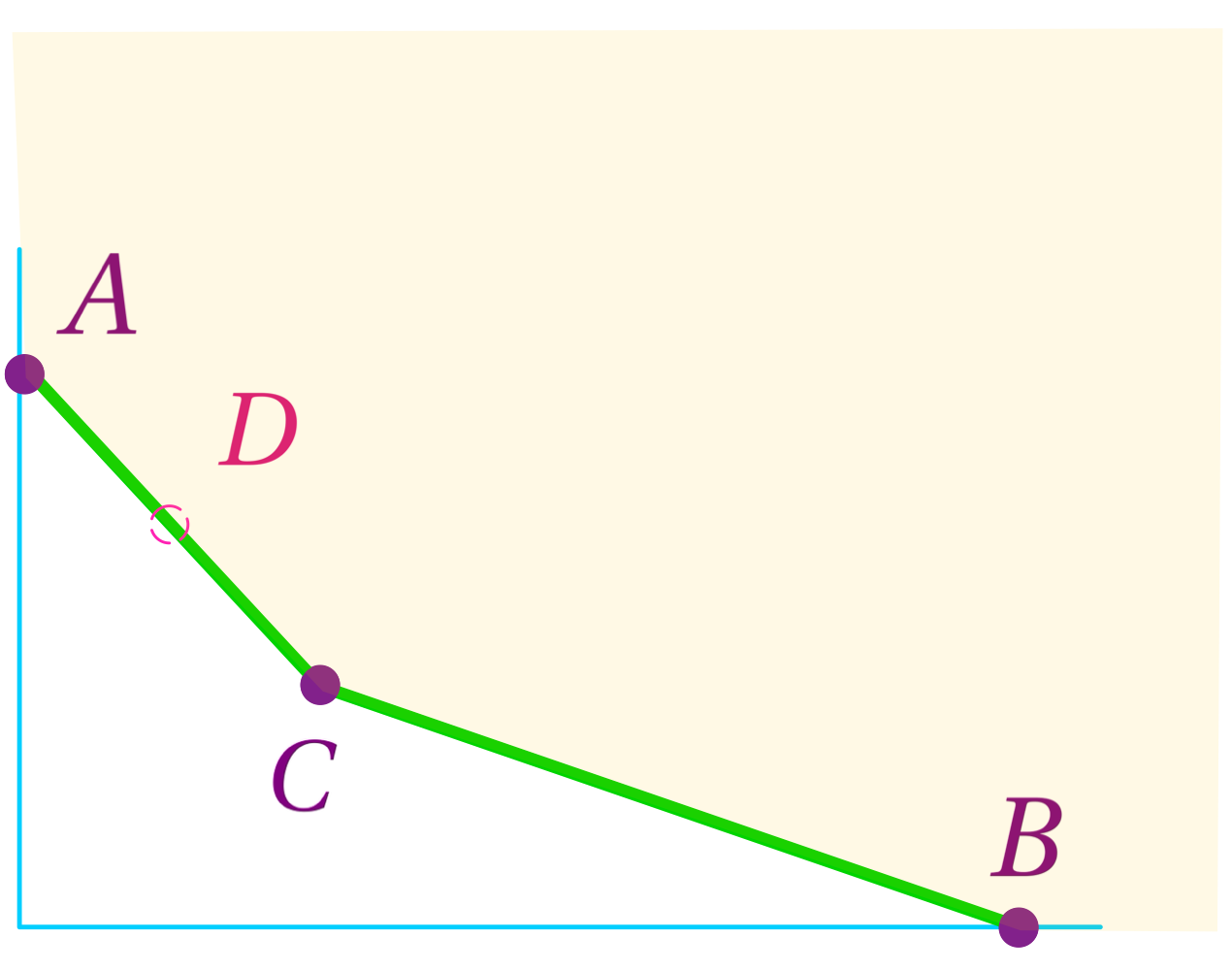}%
    }
 \caption{The minimisation of $I=(A,B,C,D,E)$ is $I_{\min}=(A,B,C)$.}
  \label{Fig:three cases}
\end{figure}

\end{proof}

As another corollary, we prove Theorem \ref{thm: normalconConvexhulldecomposition}:

\begin{theorem} \label{thm:Proof normalconConvexhulldecomposition}
  Let $I$ be an $r$-dimensional monomial ideal with singularity at the origin. The maximal dimensional irreducible components of the normal cone of $\text{\:}I$ correspond to the polytope decomposition of the lower boundary of the convex hull of $\text{\:}I$ in each $(N-j)$-coordinate plane for $0\leq j \leq r$.
\end{theorem}
\begin{proof}

From \cite[Section 3.3]{Fulton-Toric-}, we know that the rays of the associated fan to the normal cone correspond to the torus invariant exceptional divisors. Then, from  \cite[Corollary 2.53]{Telen-Toric}, these rays correspond to the maximal dimensional components of the normal cone. So, the way to construct a fan corresponding to the maximum-dimensional irreducible components of the normal cone is as follows. First, we start with the components dominating the origin, i.e., the ones embedded in the whole $N$-space and part of the lower boundary of the convex hull there, and consider their corresponding $(N-1)$-hyperplanes containing those components. Next, for the components dominating the axes we need to project the lower boundary of the convex hull onto the $(N-1)$-coordinate planes, look at the restricted convex hull there and consider the $(N-1)$-hyperplanes containing the components corresponding to the lower boundary of the convex hull there there.    We continue this process and produce $(N-1)$-hyperplanes containing the components corresponding to the lower boundary of the convex hull in coordinate planes of lower dimensions until we reach to the $(N-r)$-coordinate plane(s). Note that if $r=N-1$, then the corresponding $(N-1)$-hyperplanes  for the last step are parallel to some coordinate planes. See Figure \ref{Fig:NormalDirections} for $N=2$ and $r=N-1=1$, in which the hyperplane (line) associated to the projection onto the $x$-axis is the line $x=1$, which is parallel to the $y$-axis (as a coordinate plane). 

After constructing all the hyperplanes, the orthogonal vectors to them (which determine a fan) determine the divisors representing the maximal dimension irreducible components of the normal cone. This also aligns with \cite[Theorem 5.13]{Graffeo-Ricolfi} for $N=2$ and $r=0$.

   \begin{figure}[h]
 \subcaptionbox*{(i)}[.55\linewidth]{%
    \includegraphics[width=\linewidth]{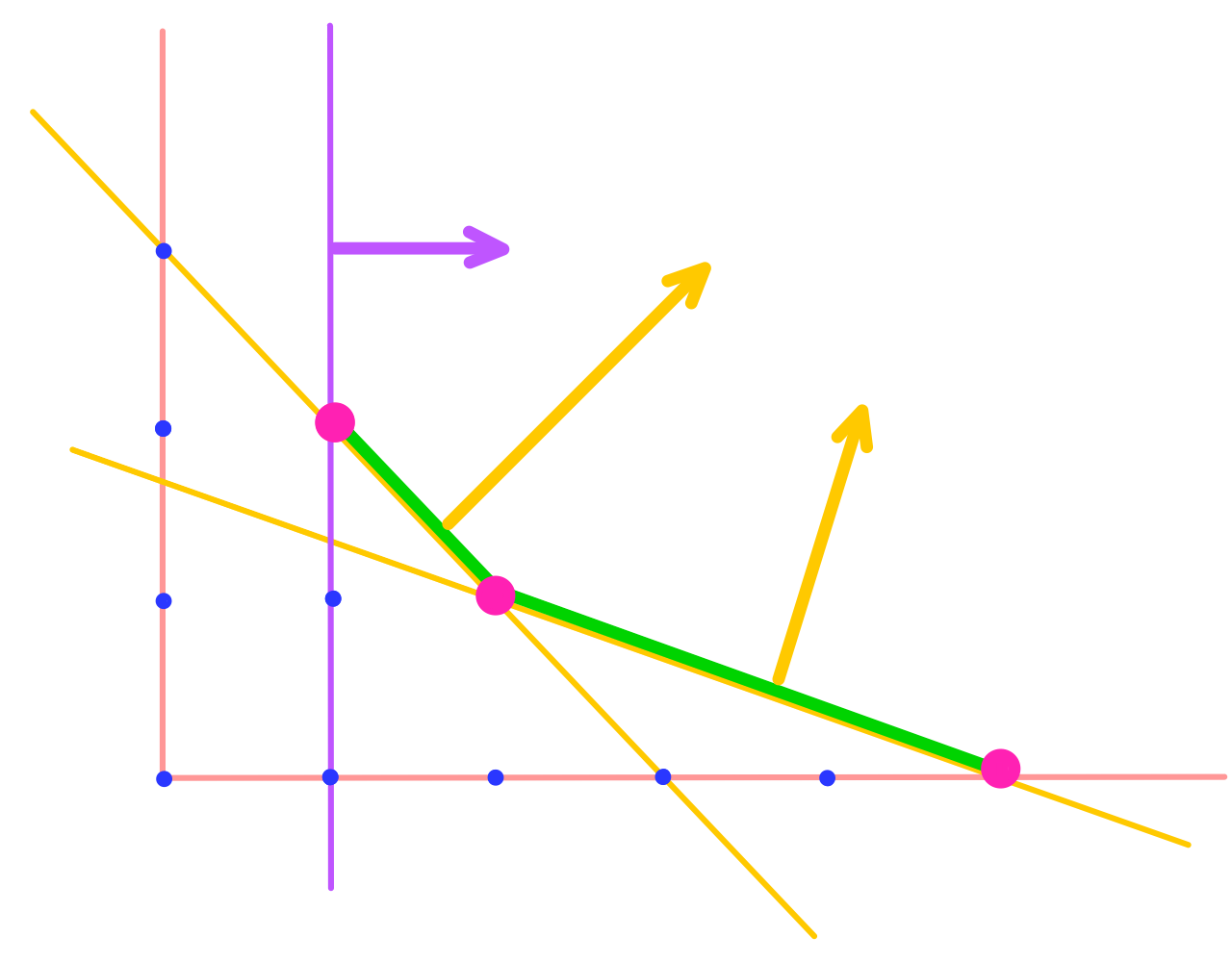}%
  }%
 \hskip12ex
  \subcaptionbox*{(ii)}[.15\linewidth]{%
    \includegraphics[width=\linewidth]{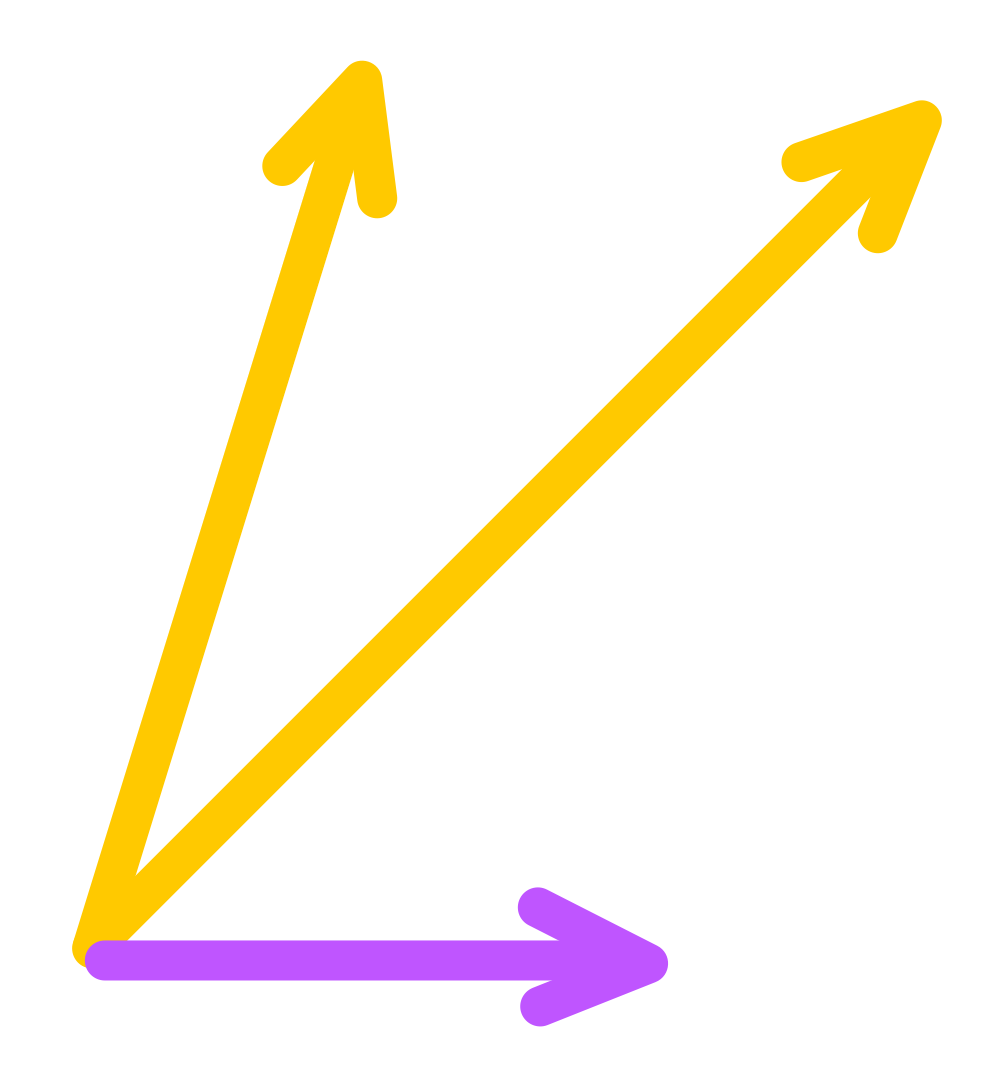}%
    
 }
  \caption{(i) The maximum-dimensional components of the normal cone associated to $I=(x^5,x^2y,xy^2)$. (ii) The fan generated by the orthogonal rays to the maximum-dimensional components. See the proof of Theorem \ref{thm:Proof normalconConvexhulldecomposition}}
  \label{Fig:NormalDirections}
\end{figure}

\end{proof}

Now, we present a  $3D$ examples for a better visualization.

\begin{example}\label{Ex:xyxzyz} Let $I=(xy,xz,yz)$ be an ideal in $\BC[x,y,z]$. In Figure \ref{Fig:xyxzyzPlanes}, apart from the hyperplane containing the triangle of the generators, there are three non-coordinate hyperplanes. Hence, the fan of the normal cone in $\BR^3$ has four rays corresponding to four maximal dimension irreducible components of the normal cone.  Also, see Example \ref{Ex: xyz}.

 \begin{figure}[h]
 \centering
    \includegraphics[width=13cm]{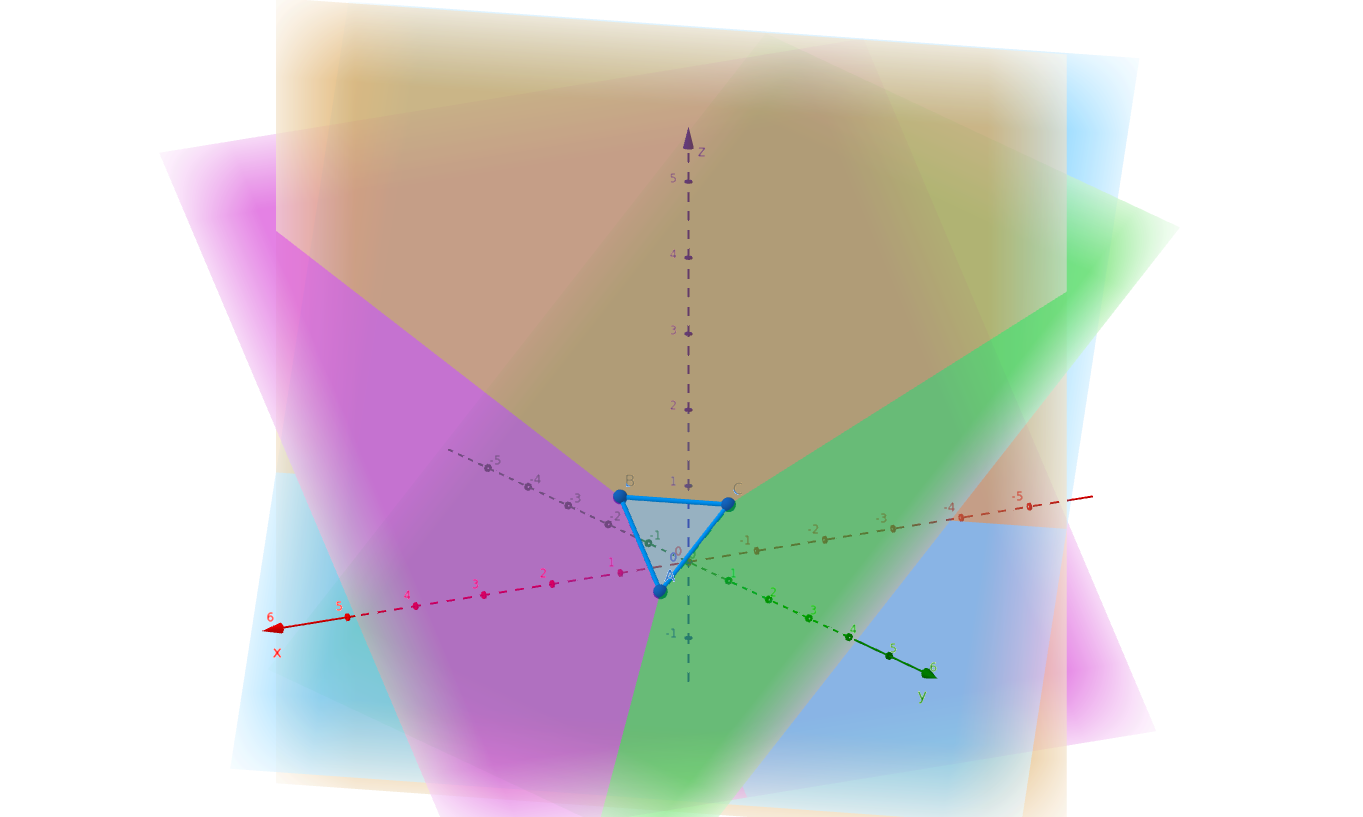}
    \caption{The lower boundary of the convex hull of the normal cone of $I=(xy,xz,yz)$ has four component corresponding to the illustrated four non-coordinate planes. See Example \ref{Ex:xyxzyz}}
  \label{Fig:xyxzyzPlanes}
\end{figure}

\end{example}

\section{A combinatorial version of the quotient formula}

 \subsection{Volume interpretation} \label{Section:VolumeInterpretation} In this sub-section, we interpret the ingredients of the quotient formula in terms of the information about the corresponding parallelepipeds.  Note that using this interpretation, we can directly compute the Behrend number without using the quotient formula. First, we have the following straightforward lemma.

 \begin{lemma}\label{primarydecomplemma}
     Let $R$ be a Noetherian ring and let $\{p_1,\ldots,p_l\}$ be minimal primes in $\mathrm{Ass}(R)$. Let $Q_i$ be a $ p_i$-primary ideal corresponding to $p_i$ with maximal dimension.
     Then, we have
     \begin{align*}
       \length(\frac{R}{Q_i})= \length(\frac{F^{-1}R}{Q_i}),
     \end{align*}
     where $F=k[y_1,\ldots,y_m]\setminus0$, where $k[y_1,\ldots,y_m]$ is the polynomial ring in the Noether normalization lemma (Lemma \ref{noetherNormalization}).

 \end{lemma}
  
  \begin{theorem}
  \label{combinatorialinterpretation}
Let $I$ be a monomial ideal of dimension $r$ in $A=\BC[x_1,\ldots,x_N]$, and let $p$ be the origin as a singular point of $X=\Spec(A/I)$. Then 
\begin{align*}
v_{X}(\mathbf{o})=&\sum_{j=1}^{n_i(I)}\sum_{i=0}^{r}(-1)^i\gcd_k(\vol(\text{parallelepipeds over  ${T^k_{i,j}}))$}\mathrm{Eu}(\big[\pi(\CN_{i,j})\big]),
\end{align*}
 where $\{T^k_{i,j}\}_k$ is a triangulation of $P_{i,j}$, where the latter is the projection of $\text{\:}\mathrm{Conv}(I)$ onto the $(N-j)$-coordinate plane for every fixed $j$. Also, $\CN_{i,j}$ are the maximal-dimensional irreducible components of $\CN(I)$, and $\pi:\CN(I)\to X$ is the canonical projection, and $\:{}\mathrm{Eu}$ is the local Euler obstruction. Also, $P_{i,j}$ corresponds to an irreducible component of the normal cone with maximal dimension, dominating a $j$-coordinate plane component of $\text{\:}\Spec(A/I)$ in $\BA^N$. 

In particular, if $I$ is zero-dimensional, we have
\begin{align*}
v_X&=\sum_{j=1}^{n(I)}v_j(I)=\sum_{j=1}^{n(I)}\gcd_k(\vol(\text{parallelepipeds over $T^k_j$}))\mathrm{Eu}(\big[\pi(\CN_{j})\big]),\\
\mult(I)&=\sum_{j=1}^{n(I)}\mult_j(I)=\sum_{j=1}^{n(I)}\vol(\text{parallelepiped over $P_j$}),\\
D_j&= (\text{number of $(N-1)$-simplices with normalised volume} \\&\text{equal to the $\gcd$ that can be accommodated in $P_j$}),
\end{align*}
where $P_j$ is a $(N-1)$-polytope corresponding to each irreducible component, and $\CN_j$ are the maximal-dimensional irreducible components of $\CN(I)$. 

\end{theorem}

\begin{proof} From the proof of Theorem \ref{v as quotien a of mult (reducedness condition)} and Theorem \ref{thm:Proof normalconConvexhulldecomposition}, the irreducible components of maximal dimension of the normal cone of $I$ are given by the maximum-dimensional components of the component $\CN_{i,j}$ corresponding to  the $p_i$-primary decomposition of the defining ideal of $\CN(I)$. To compute this using Macaulay2, see Appendix \ref{App: M2Codes}. 
These components dominate the $i$-coordinate plane(s),  which are projected onto the $(N-i)$-coordinate planes for $0 \leq i \leq r$.

 Let $r=0$. We follow the notation in the proof of Theorem \ref{v as quotien a of mult (reducedness condition)}. For volume interpretation, we correspond each monomial term $x_1^{a_{1}}x_2^{a_{2}}\ldots x_N^{a_{N}}$ in $I$ to a point $(a_{1},a_{2},\ldots,a_{N})$ in $\BR^N$ (which can be thought of as a point in $\BZ^N$). Note that for a monomial ideal $I$, the number $n(I)$ is the number of components (polytopes) in the lower boundary of the convex hull $\mathrm{Conv}(I)$.  The defining hyperplanes of the components of the lower boundary of the convex hull can be seen as orthogonal hyperplanes to the rays of associated toric fan, and each ray corresponds to an exceptional divisor of $\Bl_{\mathbf{o}}(X)$, where $\mathbf{o}$ is the origin as a singularity point of $X$. 
So, by Theorem \ref{thm:Proof normalconConvexhulldecomposition}, there is a one-to-one correspondence between the irreducible component of the normal cone corresponding to $R_j$ and the lower boundary of the irreducible components of the convex hull of $I$.

Following the notations in the proof of Theorem \ref{v as quotien a of mult (reducedness condition)} and Lemma \ref{primarydecomplemma}, we have
 \begin{align*}
 \dim_{K(R_j)}R_j=& \length(R_j)=\length(\frac{F^{-1}R}{Q_j})=\length(\frac{R}{Q_j}),\end{align*} 
  where $Q_j$ is a $p_j$-primary, so that $\{p_j\}$ are minimal primes in $\mathrm{Ass}(F^{-1}R)$.

 Let $\sigma$ be an embedding from $R$ into some $\BR^m$, where $m\geq N$. Note that $\BR^m$ contains the lattice $\BZ^m$, which in turn contains the sub-lattice $\BZ^N$, since $m\geq N$. Now, since the maximal dimension components of the normal cone correspond to ${Q_j}$, then  $\sigma$ maps $Q_j$ to a an \textit{ideal lattice} $\sigma(Q_j)=\Lambda_j$ \footnote{An \textit{ideal lattice} for a ring $R$ under an
embedding $\sigma$ into an $m$-dimensional real vector space which is a group isomorphism to $\BZ^m$, is the lattice $\sigma(I)$, where $I$  is an ideal in $R$ (see \cite[Section 1.3]{LPR})}.  Hence,  $\sigma$ is an additive isomorphism on $Q_j$. Since $R$ is $N$-dimensional,  $\Lambda_j$ is a sub-lattice of $\BZ^N$ (which, in turn, is a sub-lattice of $\BZ^m$) generated by the points in $\BZ^N$ corresponding to the lower boundary component of the convex hull of $I$ associated to $Q_j$.

Therefore, using Theorem \ref{thm:Proof normalconConvexhulldecomposition}, under $\sigma$, the restriction of $v_X$ to the $j$-th component of the normal cone, $v_{X,j}$, is
\begin{align*}
v_{X,j}=&\length (\frac{R}{Q_j})\\=&\length(\text{the maximum component of the normal cone corresponding to $Q_j$})\\=&[\BZ^N:\Lambda_j],
 \end{align*}
as $Q_j$ is $p_j$-primary, for $p_j$ minimal (isolated) prime. 

The key point here is that $[\BZ^N: \Lambda_j]$ is the $\gcd$ of all possible (normalised) volumes of $N$-simplices associated to $R_j$, considering the associated points plus the origin (by \cite[Lemma 2.4]{Blanco-Santos}; also it can be obtained from the proof of \cite[Theorem 2.2]{Stanley}). Hence, they are volumes of corresponding parallelepipeds.

For $\mult(I)$, first note that if $\Lambda$ is a sub-lattice of $\BZ^N$, and $\Gamma$ is a sub-lattice of $\Lambda$, then $[\BZ^N:\Lambda]\mid[\BZ^N:\Gamma]$  (by a straightforward argument in \cite[Lemma 2.3(1)]{Blanco-Santos}). Now, let $\Lambda_i$ be as above and $\oplus_{j=1}^{r_i}\Lambda'_{j}$ be a sub-lattice of $\Lambda_i$ generated by vertices of any fixed triangulation of $P_i$'s into $(N-1)$-simplices. By a similar argument as above, and from the proof of Theorem \ref{v as quotien a of mult (reducedness condition)}, we have
\begin{align*}
    \mult_i(I)=& \dim_{\BC(y_1,\ldots,y_N)}R_i=[\BZ^N:\oplus_{j=1}^{r_i}\Lambda'_{j}]=\sum_{j=1}^{r_i}[\BZ^N:\Lambda'_{j}]\\
    =&\sum_{j=1}^{r_i}\vol(\text{ $j$-th parallelepiped generated by each $(N-1)$-subsimplex in any}\\&\text{fixed triangulation of $P_i$}).
\end{align*}
As $\mult(I)=\sum_{i=1}^{n(I)}\mult_i(I)$, we have the result. 

Finally, for the degree of the special fiber,  we note that
\begin{align*}
     \mult_i(I)=& [\BZ^N:\oplus_{j=1}^{r_i}\Lambda'_{j}]=[\BZ^N:\Lambda_i][\Lambda_i:\oplus_{j=1}^{r_i}\Lambda'_{j}]=v_i(I)[\Lambda_i:\oplus_{j=1}^{r_i}\Lambda'_{j}].
\end{align*}
Therefore
\begin{align*}
    \sum_{j=1}^{r_i}&\vol(\text{ $j$-th parallelepiped generated by each $(N-1)$-subsimplex in any 
    triangulation of $P_i$})\\=&(\gcd(\vol(\text{parallelepipeds generated by each $N$ possible points in $P_i$})))([\Lambda_i:\oplus_{j=1}^{r_i}\Lambda'_{j}]),
\end{align*}
hence
\begin{align*}
   D_i =&[\Lambda_i:\oplus_{j=1}^{r_i}\Lambda'_{j}]=(\text{number of $(N-1)$-subsimplices with normalised volume} \\&\text{equal to the $\gcd$ that can be accommodated in $P_i$}),
\end{align*}

For $r>0$, to evaluate the Behrend function at the origin, we need to consider the alternating sum based on the dimension of the maximal dimensional irreducible component of the normal cone to obtain the formula. Using Corollary  \ref{cor: restriction} and also from the proof of  Theorem  \ref{v as quotien a of mult (reducedness condition)}, and the projection $\pi: \CN(I)\to X$, for each $0\leq i \leq r$ and fixed $j$, we know that $\CN_{i,j}$ dominates the corresponding $i$-dimensional support $\pi(\CN_{i,j})$ in $X$. For monomial $I$, this means that the support can be viewed as the projection of $\CN(I)$ onto the corresponding $(N-i)$-dimensional complement coordinate plane. Hence, the length is basically the (gcd of the) volume(s) of the parallelepiped(s). This gives us the desired formula.
\end{proof}

\begin{remark}
    One can use the existing algorithms for the convex hulls in higher dimensions (for example, see \cite{Avis-Bremner-Seidel97} and \cite{Barber-Dobkin-Huhdanpaa96}).
\end{remark}

\section{Examples}
In this section, we provide examples of different flavours and calculate the Behrend function using our main results. Note that in all the examples, we have we always have $\mathrm{Eu}=1$. So, the contribution of the local Euler obstruction will not appear in the formulas.

\subsection{Zero-dimensional examples} Let us look at some examples and directly compute the Behrend number:
\begin{example}
Let $I=(x^{10},y^{6},z^8,x^4y,y^3z,x^2z^2,xz^5)$ be a zero-dimensional ideal in $A=\BC[x,y,z]$ and let $X=\Spec(A/I)$.  From Theorem \ref{combinatorialinterpretation}, we have $v_X=20+26+24+\gcd(24,24)=94$.

\begin{figure}[h]
    \includegraphics[width=7.5cm]{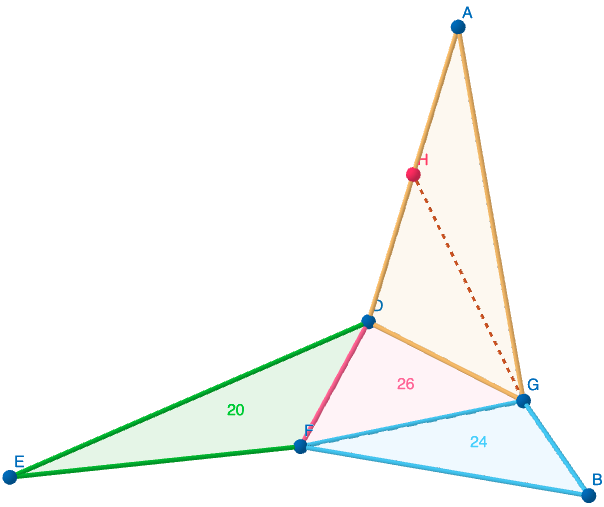}
    \caption{The monomial term $H=xz^5$ is removed in the minimisation process, but it plays a role here as it is contained in the lower boundary of the convex hull, and it gives two polytopes, each with an associated parallelepiped of volume $24$.}
    
\end{figure}

\end{example}

\begin{example}\label{Ex:114}
Let $$I=(x^{21},y^{8},z^9,x^3y^3z,xy^2z^2,xy^4z,xz^3,xy^6,x^7y^3,x^{11}y,x^9z,x^5y^4,x^7yz,x^9y^2)$$ be a zero-dimensional ideal in $A=\BC[x,y,z]$ and let $X=\Spec(A/I)$. By Theorem \ref{v as quotien a of mult (reducedness condition)}, we have \begin{align*}
    v_X=\frac{72}{1}+\frac{24}{3}+\frac{221}{17}+\frac{21}{1}=72+8+13+21=114.\end{align*}
    
    More directly, from Theorem \ref{combinatorialinterpretation}, we have

    \begin{align*}
    v_X=72+\gcd(8,8,8)+\gcd(13,13,\ldots)+21=114.\end{align*}

    See figure \ref{Fig:114}.

\begin{figure}[H]
    \includegraphics[width=7.5cm]{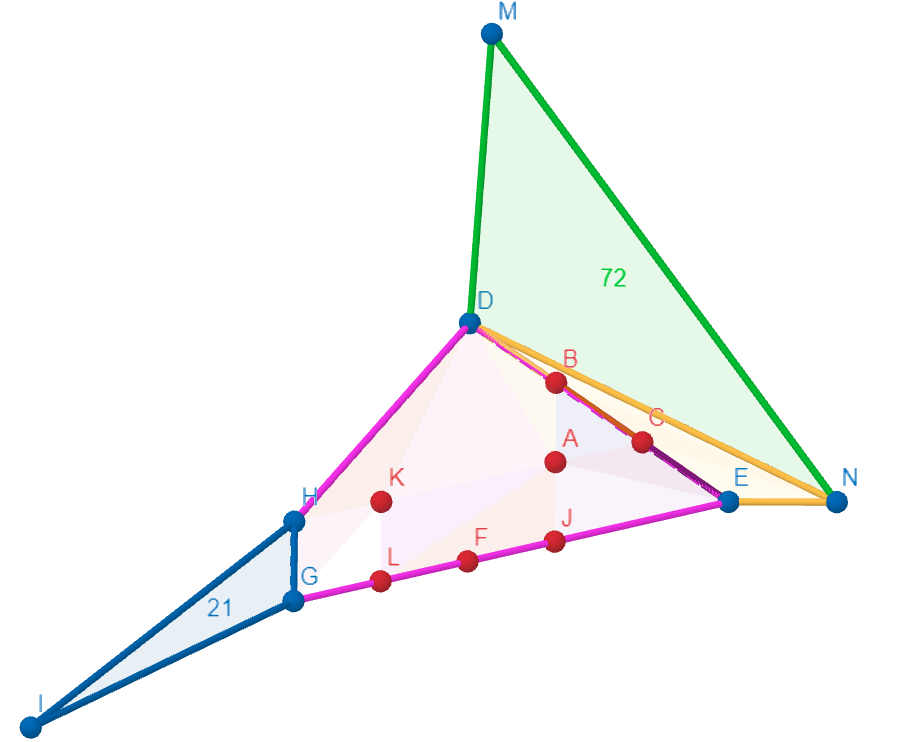}
    \caption{The monomial terms $A, B, G, K, L, F, J$ are removed in the minimisation, but they play a role in the Behrend function calculation as they are on the lower boundary of the convex hull. See Example \ref{Ex:114}.}
    \label{Fig:114}
\end{figure}

\end{example}

\begin{example}\label{Ex: Lshape} Let
\begin{align*}
    I=&(y^{17},xy^{14},x^2y^{12},x^3y^{10},x^5y^6,x^7y^4,x^9y^3,x^{11}y^2,x^{14}y,x^{16})
\end{align*}
be an ideal in $A$ and let $X=\Spec(A/I)$. Then, we have
\begin{align*}
    I_{\min}=&(y^{17},xy^{14},x^5y^6,x^7y^4,x^{11}y^2,x^{16}).
\end{align*}

As mentioned above, we can compute the Behrend number as follows. 
\begin{align*}
v_X&=\sum_{j=1}^{n(I)=5}\frac{\mult_{j}}{D_j}=\frac{17}{1}+\frac{64}{4}+\frac{22}{1}+\frac{30}{2}+\frac{32}{1}=102.\end{align*}
See Figure \ref{Fig:Lshape}.

\begin{center}
    \begin{figure}[h]
\centering
\includegraphics[width=6cm]{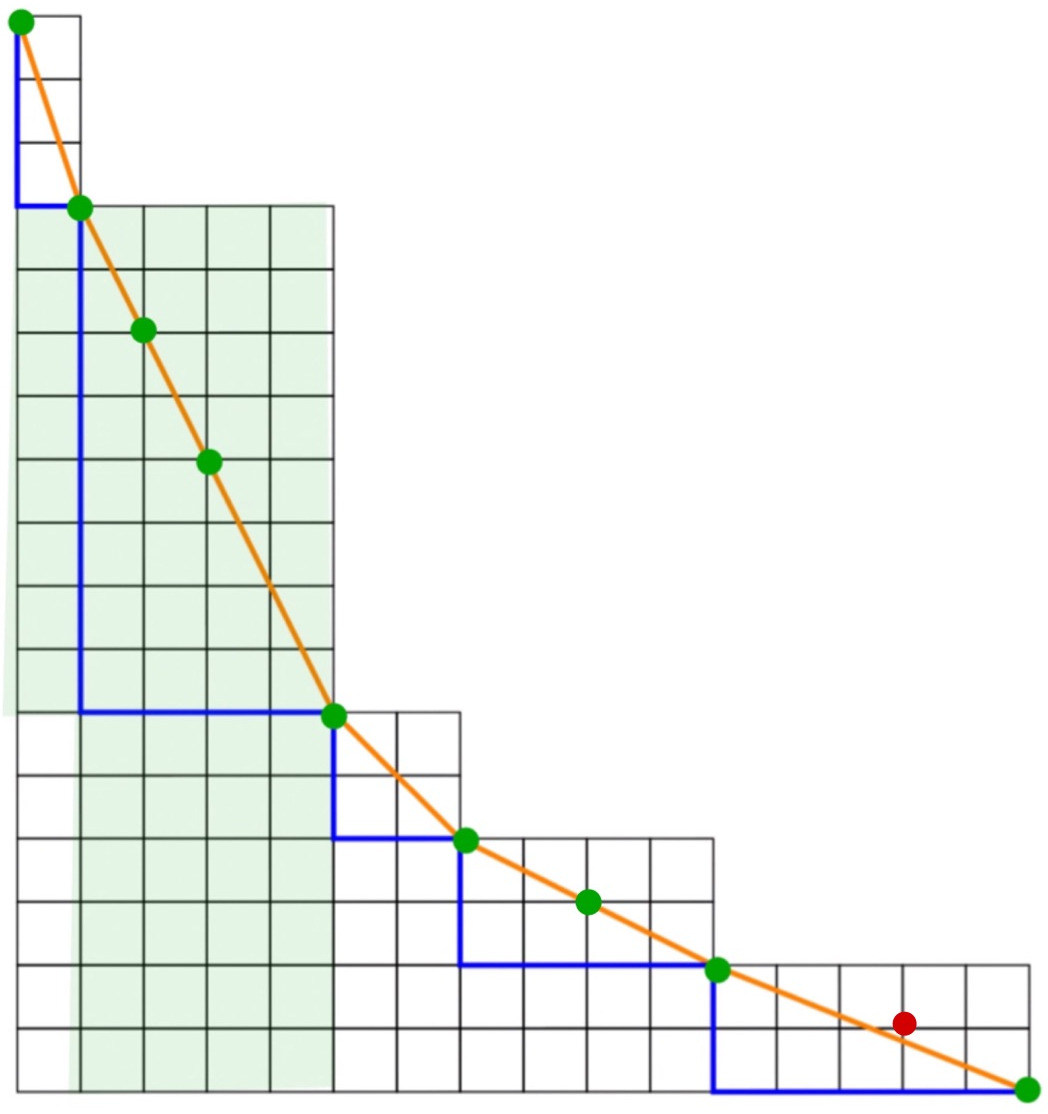}
\caption{The colored L-shape area indicates $\mult_i$ for one of the corresponding components. The length of each diagonal over the smallest possible length piece of the diagonal is $D_i$ (hence in the pictured component, we have $\mult_i=64=\text{the number of coloured squares}$, and $D_i=4$). See Example \ref{Ex: Lshape}.}
\label{Fig:Lshape}
\end{figure}
\end{center}

\end{example}

\subsection{Maximal Behrend number}
  We characterise the zero-dimensional monomial ideals with the maximum possible Behrend number (i.e., zero-dimensional ideals $I$ in $A$ with $v_X=\mult(I)$, where $X=\Spec(A/I)$):
  \begin{prop} \label{maxB}
  A zero-dimensional monomial ideal $I$ in $\BC[x_1,\ldots,x_N]$ has a maximal Behrend number if and only if the $(N-1)$-boundary of $I_{\min}$ consists only of $(N-1)$-simplices and the possible removed monomials in the minimisation process lie strictly above the $(N-1)$-boundary.
  \end{prop}
  \begin{proof}
     Assuming the conditions in the statement and from the proof of Theorem \ref{combinatorialinterpretation}, we observe that there is only one option for the $\vol$ associated with each $P_i$, which is the volume of the parallelepiped over $P_i$ itself. So, $\gcd$ is trivial and maximal, which in turn gives the maximal Behrend number.
  \end{proof}

  \begin{example}
      Let \begin{align*}
          I=(x^{11},y^{13},z^{10},x^2z^8,x^3z^4,y^7z^5,xy^6,x^2yz^4,x^5y^4z^3)
      \end{align*} be a zero-dimensional ideal in $A=\BC[x,y,z]$ and $X=\Spec(A/I)$. One can easily check that the associated $2$-polytope (or polygon) decomposition only consists of triangles, and the monomials $x^2z^8$, $y^7z^5$ and $x^5y^4z^4$, which are removed in the minimisation process, lie strictly above the lower boundary of the associated convex hull. Therefore, by Proposition \ref{maxB}, $I$ has maximal Behrend number ($v_X=\mult(I)=264+28+30+110+130=562$).  See Figure \ref{Fig:562}.
      \begin{figure}[h]
    \includegraphics[width=6.5cm]{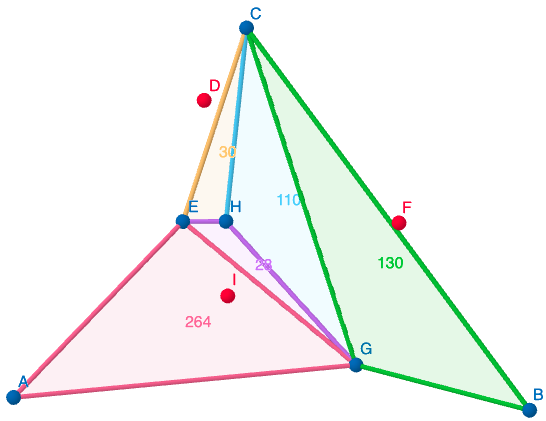}
    \caption{The red points corresponding to the removed monomials are strictly above the surface (lower boundary of the convex hull)}
    \label{Fig:562}
\end{figure}
  \end{example}

\subsection{Higher dimensional examples} We see a couple of examples of the Behrend function for the higher-dimensional monomial affine subschemes of $\BA^2$ and $\BA^3$, respectively.

\begin{example} Let $I=(xy^3,x^2y^2)$ be a one-dimensional ideal in $A=\BC[x,y]$, and $X=\Spec(A/I)$. Then, by Theorem \ref{combinatorialinterpretation}, the Behrend function at the origin $\mathbf{o}$ (as a singularity of $X$) is 
    \begin{align*}
        v_X(\mathbf{o})=&4-1-2=1.
    \end{align*}
    See Figure \ref{fig:2DEx}.

\begin{figure}[H]
 \subcaptionbox*{(1)}[.20\linewidth]{%
    \includegraphics[width=\linewidth]{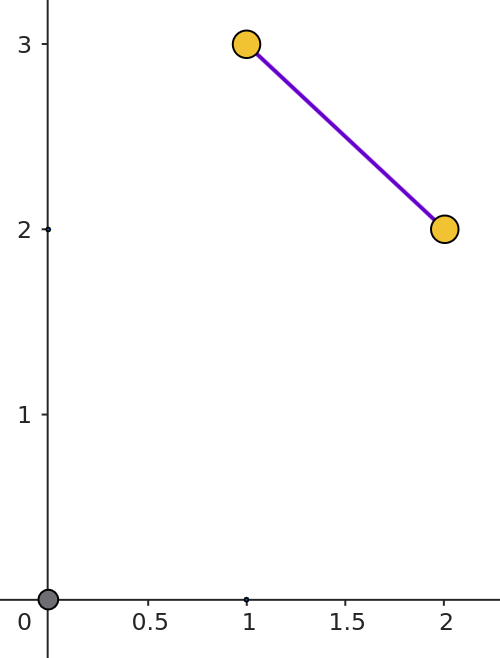}%
  }%
  \hskip10ex
  \subcaptionbox*{(2)}[.20\linewidth]{%
    \includegraphics[width=\linewidth]{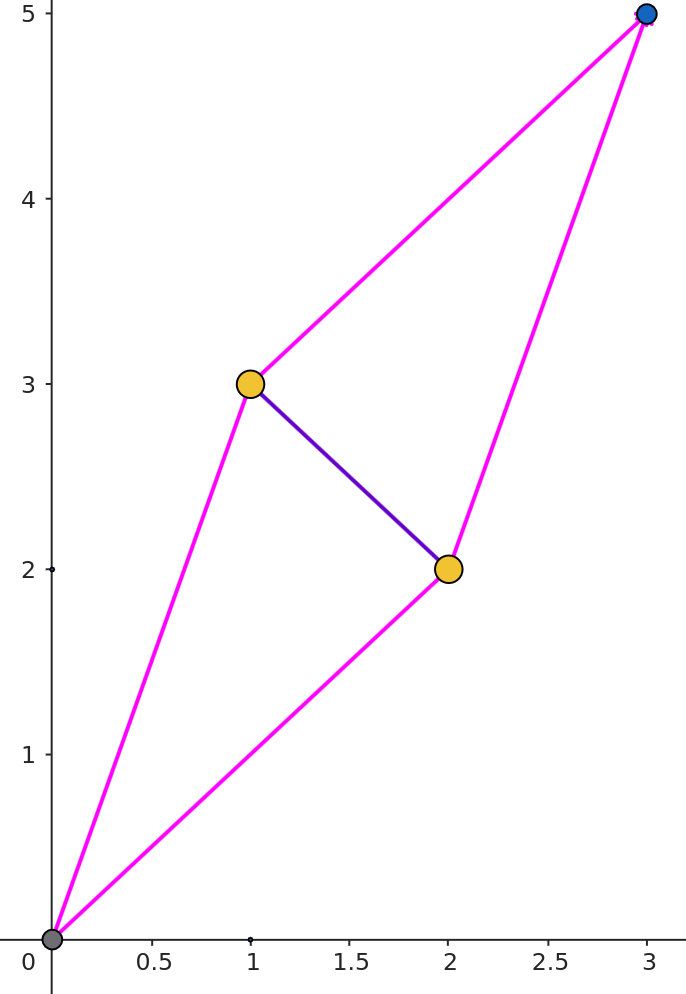}%
  }
   \hskip10ex
  \subcaptionbox*{(3)}[.20\linewidth]{%
    \includegraphics[width=\linewidth]{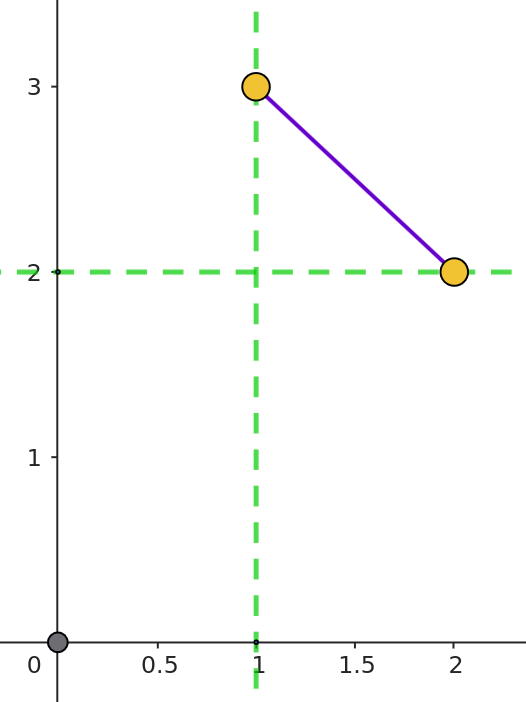}%
    }
  \caption{(1) The convex hull for $I=(xy^3,x^2y^2)$. (2) The volume of the two-dimensional parallelogram is $4$. (3) The projection on the $x$ and $y$ axes gives the lengths $1$ and $2$, respectively.}
  \label{fig:2DEx}
\end{figure}
\end{example}

\begin{example}\label{Ex: xyz}
    Let $f(x,y,z)=x^ay^bz^c$, and $X=\mathrm{Crit}(f)$. Then
    \begin{align*}
        v_X(\mathbf{o})=-1.
    \end{align*}
    This comes from Theorem \ref{combinatorialinterpretation}:
    \begin{align*}
       v_X(\mathbf{o})=& \begin{vmatrix}
a-1 & b & c\\
a & b-1 & c\\
a &b &c-1
\end{vmatrix}-\Big(\begin{vmatrix}
a & b-1 \\
a-1 & b 
\end{vmatrix}+\begin{vmatrix}
b & c-1 \\
b-1 & c 
\end{vmatrix}+\begin{vmatrix}
a & c-1 \\
a-1 & c 
\end{vmatrix}\Big)\\&+\Big((a-1)+(b-1)+(c-1)\Big)
\\=&\Big(a+b+c-1\Big)-\Big(2a+2b+2c-3\Big)+\Big((a-1)+(b-1)+(c-1)\Big)\\=&-1.
    \end{align*}
    This can be generalised to any number of variables $N$, with $f(x_1,\ldots,x_N)=x_1^{a_1}\ldots x_N^{a_N}$, for which we have $v_{\mathrm{Crit}(f)}(0)=(-1)^N$. A particular example of $N=3$ and $a=b=c=1$ is discussed in \cite[Example 2.5]{RicolfiSign}. This can also be checked by Lemma \ref{lem: CriticalLocus}.
\end{example}

\subsection{A non-monomial example} 
We provide an example of a non-monomial zero-dimensional ideal in two variables that also fits for using results in \cite{Graffeo-Ricolfi}.

\begin{example} 
    Let
    \[I=(x^4y+x^3y^2,x^5-x^3y^2,xy^5+x^3y^2+x^2y^3,x^2y^4,x^3y^3,y^7-x^3y^2-x^2y^3)\]
    be a $0$-dimensional ideal in $A=\BC[x,y]$ and $X=\Spec(A/I)$. Then, we can use the M2 codes in Appendix \ref{App: M2Codes}
to calculate the ingredients of the quotient formula and calculate the Behrend function
\[v_{X}=\frac{10}{2}+\frac{7}{1}+\frac{8}{1}+\frac{6}{1}=26.
\]
To double-check our answer, as this is a $0$-dimensional ideal in two variables for which also the conditions of \cite[Theorem 3.17 (2)]{Graffeo-Ricolfi} hold, we can write $I=K_1K_2$, where $K_1=((x+y)+(x,y))((x+y)+(x,y)^2)$ and $K_2=((x+y^2)+(x,y))((x+y^2)+(x,y)^2)((x+y^2)+(x,y)^3)$ and conclude that $v_X=5+14+2\times 2 \times 3-1\times (2+3)=26$, by \cite[Theorem 3.17 (2)]{Graffeo-Ricolfi}.
\end{example}
\subsection{Why does minimisation preserve multiplicity?}

\begin{figure}[H]
    \includegraphics[width=15cm]{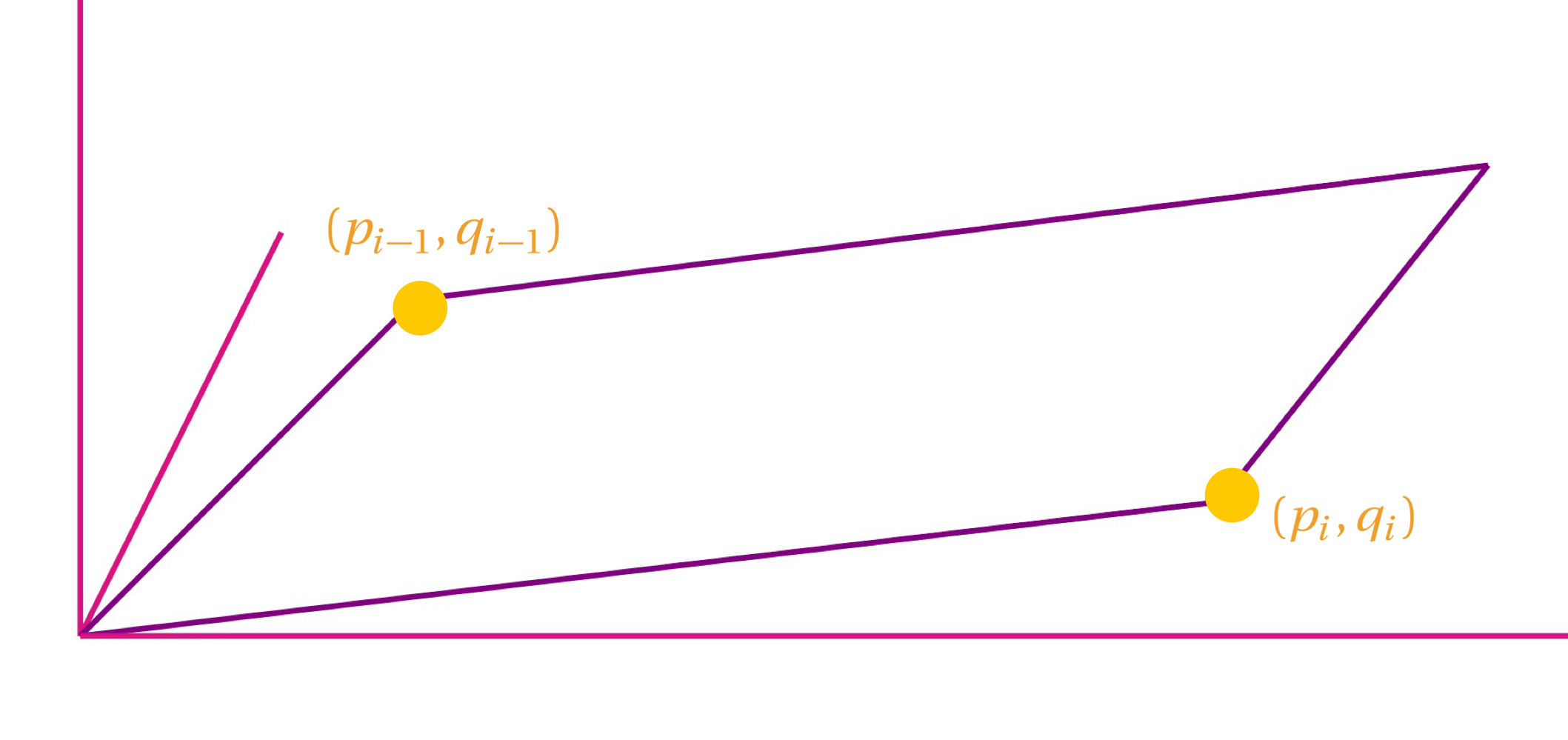}
    \caption{}
    \label{Fig:multPreserv1}
\end{figure}

As we saw in this section, the multiplicity is the sum of the volumes of parallelepipeds spanned by (sub)-simplices of each $(N-1)$-boundary associated with each monomial ideal. Without loss of generality, we answer the question for $N=2$. To answer the question, we need to check it for each $\mult_i$, i.e., for the parallelogram generated by vectors $(p_{i-1},q_{i-1})$ and $(p_i,q_i)$, for each $i$. For the sake of simplicity, we restrict to two dimensions, and the argument could be easily generalised to higher dimensions. Let $I$ be an ideal in $\BC[x,y]$, so that the parallelogram produced by two of the successive generators $x^{p_{i-1}}y^{q_{i-1}},x^{p_{i}}y^{q_{i}}$ of $I_{\min}$ is as in Figure \ref{Fig:multPreserv1}.

\begin{figure}[H]
    \includegraphics[width=15cm]{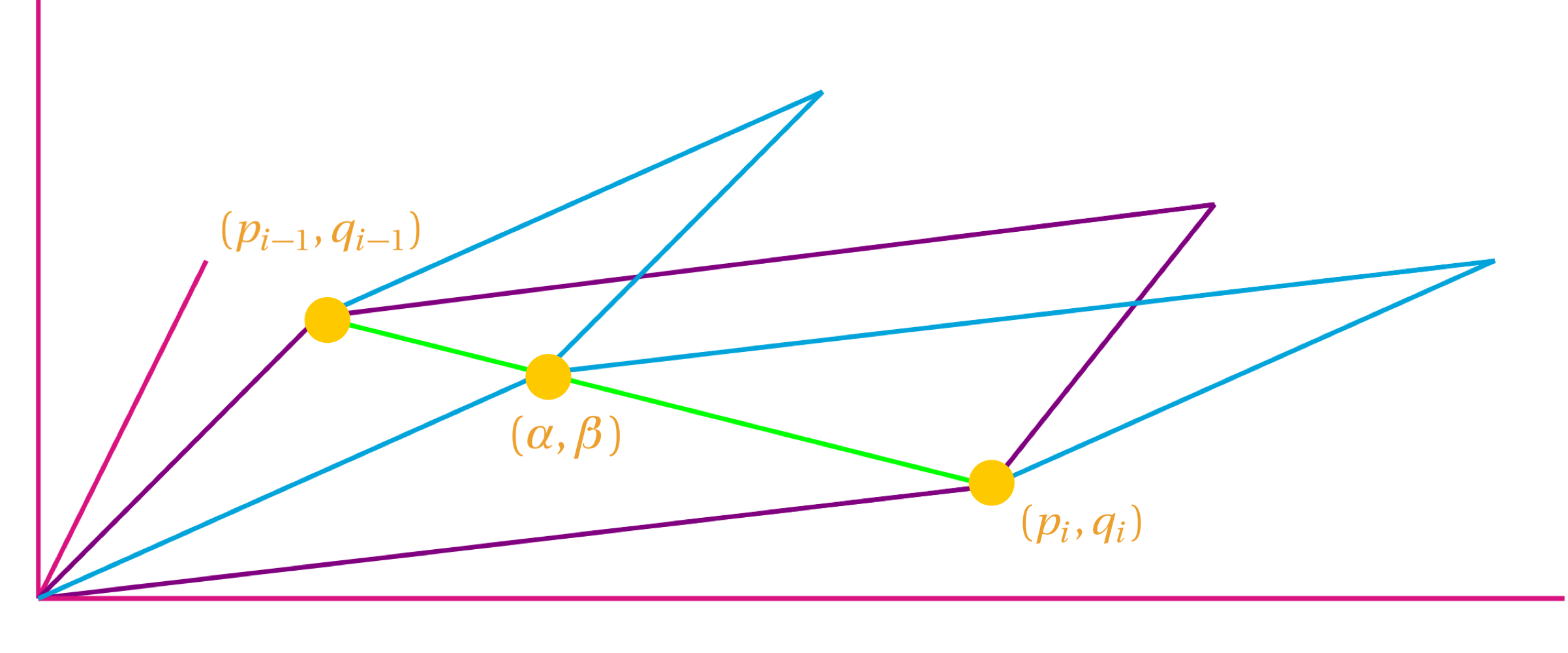}
    \caption{}
    \label{Fig:multPreserv2}
\end{figure}

  Now, assume that there is a removed generator $x^{\alpha}y^{\beta}$ on the smaller diagonal of the parallelogram, which is removed in the process of minimisation. Associated to such a point, there are two smaller parallelograms as in Figure \ref{Fig:multPreserv2}.

  We want to show that the area of the parallelogram generated by the vectors $(p_{i-1},q_{i-1})$ and $(p_i,q_i)$ is equal to the sum of the areas of the parallelograms generated by $(p_{i-1},q_{i-1})$ and $(\alpha,\beta)$, on the one hand, and $(\alpha,\beta)$ and $(p_i,q_i)$, on the other hand. We denote the area of the parallelogram generated by vectors $ v$ and $ w$ by $A(v,w)$. First note that the line from the point $(p_{i-1},q_{i-1})$ to the point $(p_{i},q_{i})$ is given by the following equation in the $xy$ plane:
  
  \begin{align*}
      y-q_{i-1}=\frac{q_{i}-q_{i-1}}{p_{i}-p_{i-1}}(x-p_{i-1}).
  \end{align*}
  As the point $(\alpha,\beta)$ is on this line, we have
  \begin{align*}
       \alpha=\frac{\beta(p_{i}-p_{i-1})-q_{i-1}(p_{i}-p_{i-1})+p_{i-1}(q_{i}-q_{i-1}))}{q_{i}-q_{i-1}}.
  \end{align*}
 Therefore, we have
  \begin{align*}
  A((p_{i-1},q_{i-1}),(\alpha,\beta))+A((\alpha,\beta),(p_i,q_i))=&(\alpha q_{i-1}-p_{i-1}\beta)+(p_{i}\beta-\alpha q_{i})\\
  =&\alpha(q_{i-1}-q_{i})+\beta(p_{i}-p_{i-1}) \\
  =&-\beta(p_{i}-p_{i-1})+q_{i-1}(p_{i}-p_{i-1})\\-&p_{i-1}(q_{i}-q_{i-1}))+\beta(p_{i}-p_{i-1})\\
  =&q_{i-1}p_i-q_ip_{i-1}\\
      =& A((p_{i-1},q_{i-1}),(p_i,q_i)),
  \end{align*}
as desired.

\appendix

\section{A generalised Pick's Theorem}
The argument in this section can be written without any mention of the Behrend function. However, as we use a similar combinatorial spirit as of the Behrend function for monomial ideals, we put the result here.
\subsection{Behrend function in terms of the equations of the lower boundary of the convex hull in three dimensions} We have the following lemma to recover the contribution of the Behrend function on each component from the equation of the plane containing the polygon.

\begin{theorem}\label{v=|d|}
Let $I$ be a monomial ideal in $A=\BC[x_1,x_2,\ldots,x_N]$ that is singular at the origin, and let $X=\Spec(A/I)$. Then, for any polygon $P_{i,j}$ in the lower boundary of the convex hull of $I$,  
we have
\begin{align*}
    v_{X,i,j}(\mathbf{o})=&|d|,
\end{align*}
where $v_{X,i,j}$ is the contribution of the Behrend function evaluated at the origin (considering with a suitable sign), and $d$ is the constant term of the primitive equation of the plane defined by $P_{i,j}$.
\end{theorem}
\begin{proof}
Let $a_1x_1+a_2x_2+\ldots+a_Nx_N+d=0$ be the primitive\footnote{$\gcd(a_1,\ldots,a_N,d)=1$} equation of the plane defined by $P_{i,j}$. We only need to restrict to the case when the polygon $P_{i,j}$ is just a ($N-1$)-simplex, and the only lattice points on it are just the vertices. From Theorem \ref{combinatorialinterpretation}, we have
$$v_{X,i,j}(\mathbf{o})=\vol(\text{parallelepiped spanned by the vertices of $P_{i,j}$}).$$

If $A_1, A_2,\ldots, A_N$ are the vertices
 (as in the Figure \ref{Fig:ABC} for $N=3$ and setting $A_1=A,A_2=B$ and $A_3=C$), then we have 
 \begin{align*}\label{area}
     v_{X,i,j}(\mathbf{o})=&\vol(\text{parallelepiped spanned by $A_1, A_2,\ldots, A_N$})=N!\vol(\text{$N$-simplex spanned by $A_1, A_2,$}\\&\text{$\ldots, A_N$ and the origin})\\=&N!(\frac{1}{N}\vol(P_{i,j})h)=(N-1)!\vol(P_{i,j})h,
 \end{align*}
 where $h$ is the distance of the origin from $P_{i,j}$, which is equal to $\frac{|d|}{\sqrt{a_1^2+\ldots+a_N^2}}$. On the other  hand, $\vol(P_{i,j})=\frac{1}{(N-1)!}|\overrightarrow{n}|$, where $\overrightarrow{n}$ is the normal vector to $P_{i,j}$. It is straightforward to check that by our assumption, we have $|\overrightarrow{n}|=\sqrt{a_1^2+a_2^2+\ldots+a_N^2}$, in this case. Hence, we have
 \begin{align*}
     v_{X,i,j}=&(N-1)!\vol(P_{i,j})h=(N-1)!(\frac{1}{(N-1)!}|\overrightarrow{n}|)(\frac{|d|}{\sqrt{a_1^2+a_2^2+\ldots+a_N^2}})
     =|d|.
 \end{align*}
 \begin{figure}[H]
 \centering
    \includegraphics[width=5cm]{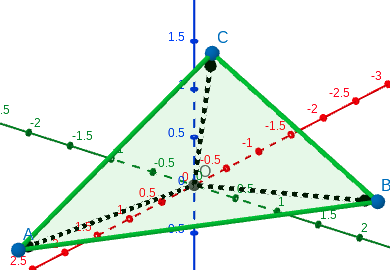}
    \caption{}
\label{Fig:ABC}
\end{figure}
    
\end{proof}

\begin{example}
    Let $I=(x^3,y^3,z^2,xy,xz,yz)$ be a zero-dimensional ideal in $\BC[x,y,z]$. Using Theorem \ref{v=|d|}, we can see that the defining equation of the three planes containing the irreducible component of the normal cone of $I$ are $x+y+z-2=0$, $2x+y+2z-3=0$ and $2x+y+2z-3=0$. Hence, the Behrend number is $2+3+3=8$. See Figure \ref{Fig:GRexample}. Note that this is the example in \cite[Section 7]{Graffeo-Ricolfi} where the authors introduce to show the difficulties in three dimensions.
        \begin{figure}[H]
 \subcaptionbox*{(1)}[.32\linewidth]{%
    \includegraphics[width=\linewidth]{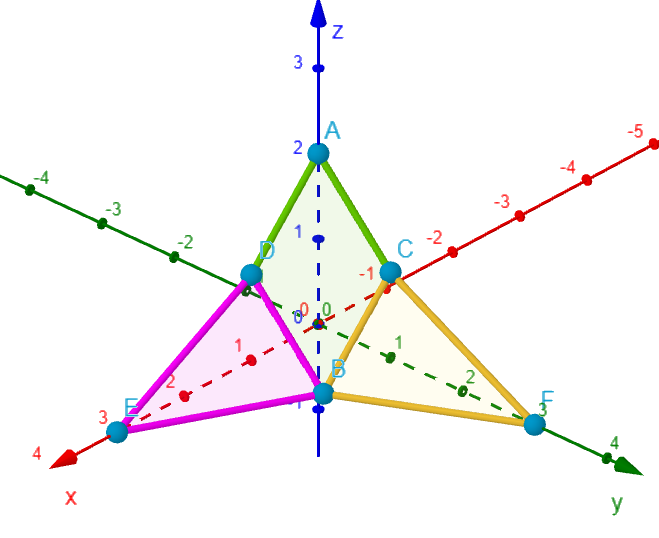}%
}%
  \hskip2ex
  \subcaptionbox*{(2)}[.62\linewidth]{%
    \includegraphics[width=\linewidth]{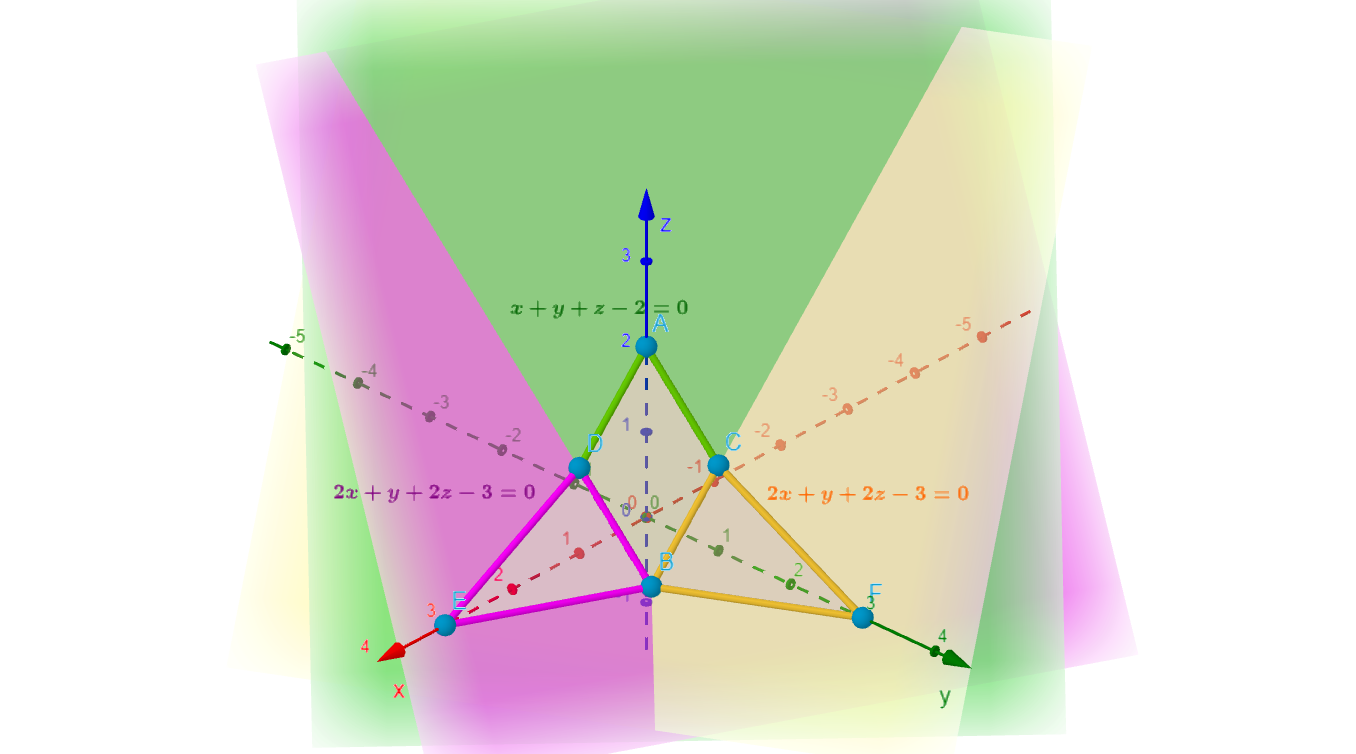}%
  }
  \caption{(1) The irreducible components of the normal cone of $I=(x^3,y^3,z^2,xy,xz,yz)$. (2) The corresponding defining planes.}
  \label{Fig:GRexample}
\end{figure}
\end{example}
\subsection{Pick's Theorem for embedded polygons} We prove a generalised version of Pick's theorem for embedded polygons in $\BR^3$. The classical Pick's theorem states the area of a plane polygon in terms of the number of the lattice points inside and on the boundary of the polygon. First, we need the following definition and lemma.

\begin{definition}[Fine triangulation]\label{Def: fineTriangulation} A \emph{fine triangulation} for a polygon, is an integral triangulation (i.e., the vertices of all the triangles are lattice points) with the maximum number of triangles.
\end{definition}

\begin{definition}[Empty simplex/empty triangle] An \emph{empty $N$-simplex} is a symplex with lattice points only on the vertices. An empty $2$-simplex is an \emph{empty triangle}.
    
\end{definition}
\begin{lemma}\label{number B I}
    Let $P$ be a polygon in the primitive plane $ax+by+cz+d=0$, and let $I(P), B(P)$ be the number of internal lattice points and the number of lattice points on the boundary of $P$, respectively. Then
    \begin{align*}
        \text{the number of empty triangles in a fine integral triangulation of $P$}=&B(P)+2I(P)-2.
    \end{align*}
    
\end{lemma}
\begin{proof}
 Without loss of generality, and up to transformations, we can assume that $c=1,a,b,d=0$. We proceed by induction on the number of lattice points which are not vertices. The base case is when $P$ a empty triangle, for which we have $I(P)=0, B(P)=3$, hence $B(P)+2I(P)-2=1$, which is obviously the number of empty triangles in $P$. Now, let $P$ be a polygon with at least one non-vertex lattice point, and let $P'$ be a polygon purely inside $P$ with $B'(P')$ lattice points on its boundary and $I'(P')$ internal lattice points, and we can assume that the statement holds for $P'$ and we want to show that $B(P)+2I(P)-2$ is also the number of empty triangles in $P$.  See Figure \ref{Fig:IBI'B'}(1).
 
 If there is no lattice point between $P$ and $P'$, then there are $B(P)+B'(P')$ empty triangles. between $P$ and $P'$: set $m_1=\min\{B(P),B'(P')\}$ and $m_2=\max\{B(P),B'(P')\}$. Then there are $m_1$ (potentially, degenerate or concave) quadrilaterals with two sides of it being a side of $P$ and a side of $P'$. This gives $2m_1$ empty triangles.   On the other hand, for each remaining $m_2-m_1$ side, there exist a unique triangle. Totally this gives $2m_1+(m_2-m_1)=m_1+m_2=B(P)+B'(P')$ empty triangles between $P$ and $P'$. 
 Now, note that each lattice points $p$ between $P$ and $P'$ adds two more empty triangles to $B(P)+B'(P')$ (see Figure \ref{Fig:IBI'B'}(2)); i.e., there are totally $B(P)+B'(P')+2(I(P)-I'(P')-B'(P'))$ empty triangles between $P$ and $P'$. 
 Combining this with the induction hypothesis (the number of empty triangles in $P'$ is $B'(P')+2I'(P')-2$), we conclude that there are $B'(P')+2I'(P')-2+B(P)+B'(P')+2(I(P)-I'(P')-B'(P'))=2I(P)+B(P)-2=$ empty triangles in $P$, as claimed.

 Note that the result holds for every plane in $\BR^3$, since linear transformation of planes in $\BR^3$ preserves the lattice points and empty triangles.

     \begin{figure}[H]
 \subcaptionbox*{(1)}[.4\linewidth]{%
    \includegraphics[width=\linewidth]{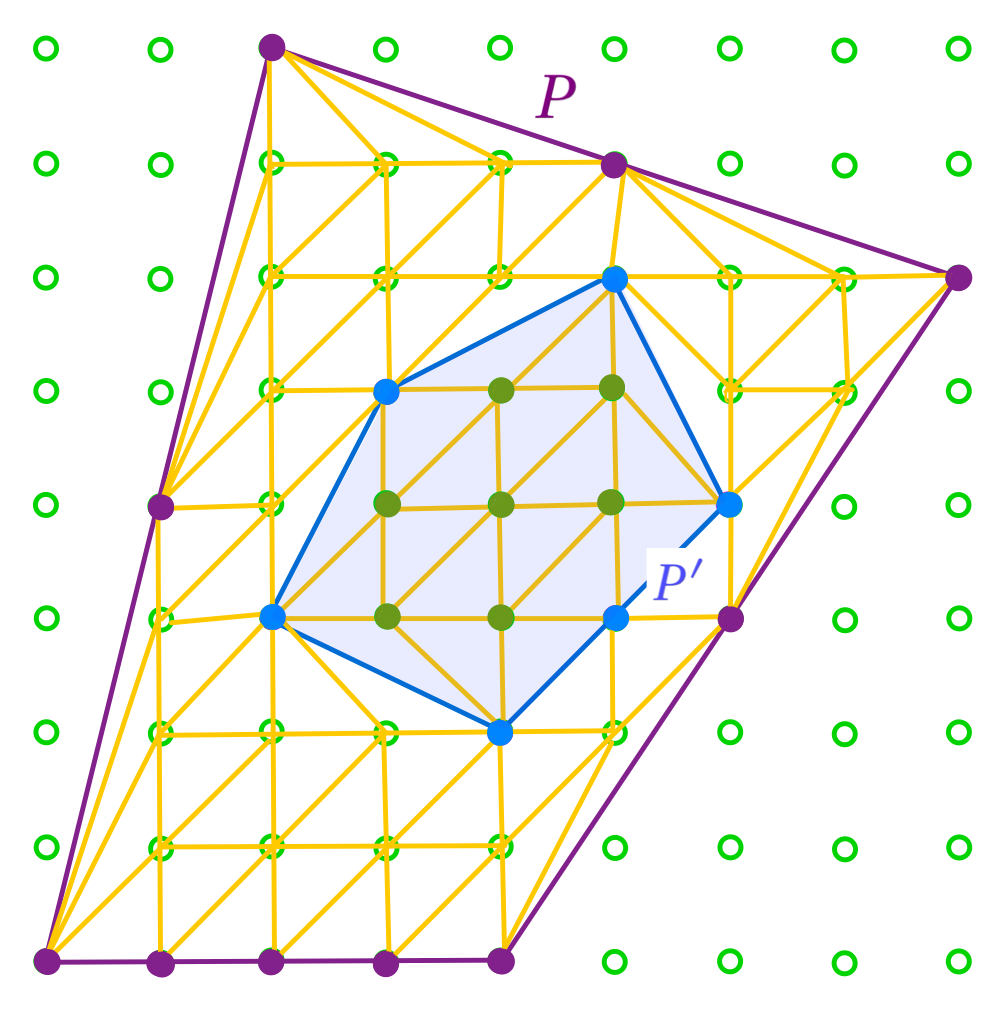}%
}%
  \hskip10ex
  \subcaptionbox*{(2)}[.38\linewidth]{%
    \includegraphics[width=\linewidth]{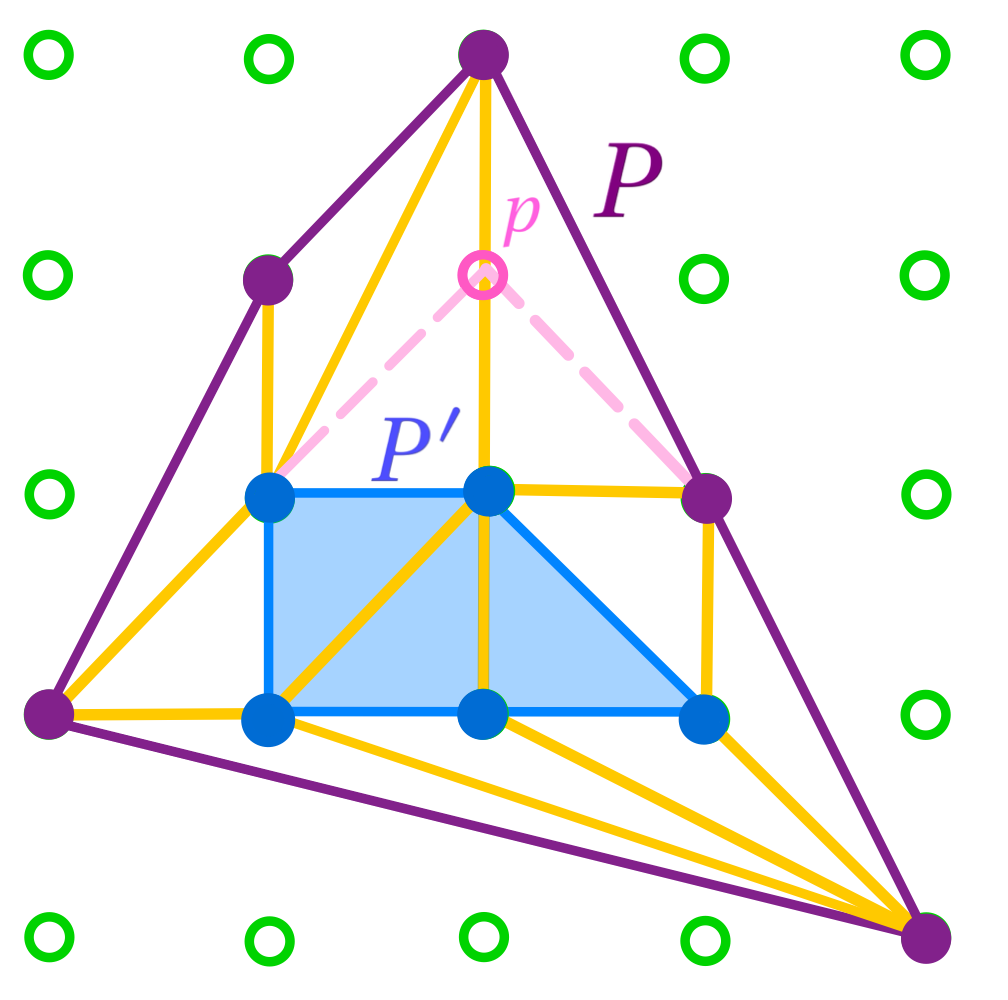}%
  }
  \caption{(1) Proof by induction. (2) Any lattice point $p$ between $P$ and $P'$ adds two additional empty triangles.}
  \label{Fig:IBI'B'}
\end{figure}

\end{proof}

\begin{example}Figure \ref{Fig:IBI'B'} shows two examples related to Lemma \ref{number B I}, for $B(P)=10, I(P)=34$ (left) and $B(P)=5, I(P)=6$ (right).

\end{example}
Now, we can prove a generalised version of Pick's Theorem:
\begin{theorem} \label{Pick's}
    Let $P$ be a plane convex polygon with some non-lattice vertices so that it can be integrally and generically embedded in $\BR^3$, i.e., there is an embedding of $P$ into $\BR^3$ with all vertices lattice points so that it is not perpendicular to any of the three coordinate planes.  Let $\widetilde{I(P)}$ and $\widetilde{B(P)}$ be the number of internal lattice points and the number of lattice points on the boundary of the embedded polygon in $\BR^3$, respectively.  Then
    \begin{align*}
        \area(P)=&\left(\frac{\widetilde{B(P)}}{2}+\widetilde{I(P)}-1\right)\left(\sqrt{a^2+b^2+c^2}\right),
    \end{align*}
    where $ax+by+cz+d=0$ is the primitive equation of the plane defined by the embedded polygon in $\BR^3$, and $(a,b,c)$ is primitive, and $d\neq 0$.
\end{theorem}

\begin{proof}
From the proof of Theorem \ref{v=|d|} for $N=3$, the volume of each empty triangle in $P$ with primitive equation $ax+by+cz+d=0$ is $\sqrt{a^2+b^2+c^2}/2$, and from  Lemma \ref{number B I}, the number of such empty triangles in $P$ is $\widetilde{B(P)}+2\widetilde{I(P)}-2$. Combining these implies the claim.

\end{proof}

\section{Macaulay2 codes}\label{App: M2Codes}
The codes in this appendix are due to Mike Stillman in April 2023. The codes compute the ingredients of the quotient formula for a given ideal, which defines our affine scheme. Here, we consider the example of the one-dimensional subscheme of $\BA^3$ defined by the ideal $I=(xy,yz,xz)$. Here, \verb+L+ is the same as in the proof of Theorem \ref{v as quotien a of mult (reducedness condition)}, and \verb+Qs+ is the same as $Q_{i,j}$ or $Q_j$, and \verb+Ps+ is the same as $p_{i,j}$ or $p_j$ (as in the proof of Theorem \ref{v as quotien a of mult (reducedness condition)}).
\vspace{0.2cm}

\begin{verbatim}
R=QQ[x,y,z];
I=ideal (x*y,y*z,x*z);
J=trim I;
R1 = normalCone J;
S = first flattenRing R1;
Rsg = (coefficientRing S)[gens S];
L = sub(trim ideal S, Rsg);
Qs = primaryDecomposition L;
maxdim = max (Qs/dim);
Qs = select(Qs, q -> dim q == maxdim);
Ps = Qs/radical;
Qs/degree
Ps/degree
\end{verbatim}

\bibliographystyle{amsplain-nodash} 
\bibliography{bib}

\end{document}